\documentclass{article}
\usepackage{amssymb,amsmath}
\usepackage{amsfonts}
\usepackage{latexsym}
\newtheorem{Lem}{Lemma}
\newtheorem{Thm}{Theorem}
\newtheorem{Pro}{Proposition}
\newtheorem{Rem}{Remark}
\newtheorem{Cor}{Corollary}
\newtheorem{Exa}{Example}

\def\C{{\mathbb C}}
\def\R{{\mathbb R}}

\def\M{{\mathbb M}} 

\def\Hh{{\mathbb H}}
\def\M{{\mathbb M}}
\def\K{{\mathbb K}} 
\def\A{{\mathcal A}}
\def\B{{\mathcal B}}
\def\x{{\bf x}}
\def\y{{\bf y}}
 
\def\s{{\bf  s}} 
\def\q{{\bf q}}
\def\H{{\mathfrak{H}}}
 
\title{NORMAL OPERATORS IN REAL AND QUATERNIONIC HILBERT SPACES}

\author{Florian-Horia Vasilescu\\
\small Department of Mathematics, University of Lille,\\
\small 59655 Villeneuve d'Ascq, France\\
\small e-mail: florian.vasilescu@univ-lille.fr}
\date{}

\begin{document}

\maketitle

{\it Keywords: } spectral measures; complex extensions of real operators; spectral representations; real and quaternionic normal 
operators

{\it AMS Subject Classification 2020:} 47B15; 47A60; 47S05

\begin{abstract} 

A new approach to normal operators in real Hilbert spaces is discussed,
and a spectral representation is obtained, derived directly from the 
complex case. The results are then applied to quaternionic normal operators, regarded as a special class of real normal operators. This
point of view  allows us to consider their spectrum and associated  measures to be defined  on subsets of the complex plane, in a classical
manner.
\end{abstract}

\section{Introduction}

The concept of normal operator, bounded or not, has a central role
in operator theory and its applications, and it has been preponderantly used in the framework of complex Hilbert spaces. Nevertheless, normal operators in real Hilbert spaces have been studied by several authors,
including Goodrich \cite{Goo}, Agrawal and  Kulkarni \cite{AgKu}, 
Viswanath \cite{Vis0}, and Oreshina \cite{Ore}, to quote some of them,  but the oldest contribution in this respect seems to be that of Wong in \cite{Won}. Let us briefly refer to the contributions of these authors. 

In the paper \cite{Goo}, Goodrich exploits the complexification of a 
real Hilbert space, to obtain a spectral representation  of a real 
normal operator. In fact, he expresses  a real normal operator as a sum of two integrals, via a concept called
spectral pair, which consists of finite variation and countably additive measures  in the strong operator topology, combined with some 
trigonometric functions. This result is then used to  show that every real normal operator  is orthogonally equivalent to an orthogonal sum of operators acting on spaces of square integrable functions. 

Agrawal and Kulkarni (see \cite{AgKu})  prove a spectral theorem for a given normal operator on a real Hilbert space, similar to the complex case, by using the techniques of real Banach algebras. Specifically, they use  an analogue of the 
Gelfand-Naimark Theorem by Arens-Kaplansky (see \cite{KulLim}), and a real version of the Riesz Representation Theorem by Grzesiak \cite{Grz}. 

As mentioned in \cite{AgKu}, a version of the spectral theorem for a normal operator acting on a real Hilbert spaces was proved by Viswanath using some concepts from the measure theory. The paper \cite{Vis0} has been out of reach of the present author.  

Oreshina (see \cite{Ore}) uses the representation of a
normal operator as a sum of a self-adjoint and a skew-adjoint operator 
obtaining  the spectral decomposition as the sum of two integrals.  Because this approach fails the construction of the functional calculus, a second version of the spectral theorem is proved to recapture such a construction.

In this work, we basically show  that a very consistent information concerning a real normal
operator can be directly obtained by using the corresponding properties of its complex extension, which is also a normal one.  The assumed properties 
are well known and unanimously accepted, a valuable source being  the functional analysis monograph by Rudin  \cite{Rud} (see also \cite{DuSc} and \cite{Sch}).  Our approach to
real case not only brings some new information in a simplifid manner but it is also indispensable for the further development in the quaternionic context, where the quaternionic normal operators are regarded as a special class of real normal operators.

In our discussion concerning the real normal operators, we start with the bounded case, the unbounded one being subsequently approached, using some arguments from the former case. Our main results in this respect are Theorem \ref{bdfc0} for the bounded case, and Theorem \ref{ubdfc} for the unbounded one. As in the paper \cite{Goo}, we exploit the existence of a spectral measure for the complex normal extension of the given real normal operator but our arguments are totally different. We show that the restriction of the complex 
construction to the original real Hilbert space leads to the desired properties. Our Theorem \ref{bdfc0} proves the existence of what may be called a {\it spectral representation}, which is in fact a functional calculus with a large classes of Borel functions associated to the
given real normal operator. This representation behaves like a 
spectral measure but it is  a slightly more general concept, however
having a unique family of associated scalar measures (see Corolary \ref{bdfc1}). A natural associated spectral measure is not strong enough to be applied to all
involved functions, as shown by Example \ref{con_exa}. Nevertheless, 
using some elements  of local spectral theory, 
specifically the uniqueness of spectral capacities (see \cite{Vas0}), under some general conditions
we prove uniqueness results concerning the real spectral measure attached to a real normal
operator (see Propositions \ref{uniqueness} and \ref{unb_uniqueness}).

The case of unbounded real normal operators stated by Theorem \ref{ubdfc} is obtained using partially  our Theorem 1, as well as other results from the monograph \cite{Rud}.

The first investigations in the framework of quaternionic Hilbert spaces seemingly go back to  \cite{Tei} (see also \cite{Vis}, \cite{AlCoKi}, \cite{CoGaKi} etc.). 

In the paper \cite{Vis}, the author obtains a spectral theorem, a
functional calculus and a multiplicity theory for normal operators acting in quaternionic Hilbert spaces by refining the corresponding complex methods. The algebra of all operators and unitary representations of locally compact abeliang groups in quaternionic Hilbert space are also studied. 

The authors of the paper \cite{AlCoKi} obtain a spectral theorem in 
quaternionic Hilbert spaces using the notion called $S$-spectrum,
which is a subset in the quaternionic algebra $\Hh$. A similar point of view is adopted by the authors of the monograph \cite{CoGaKi}. 

In the paper \cite{GhMoPe1}, new mathematical objects are introduced and studied. They are called intertwining quaternionic projection-valued measures, or shortly iqPVMs.

The authors of the paper \cite{RaKu} deal with the position operator
from quantum mechanics, in the quaternionic setting, which is normal, proving a version of the multiplication form  of the spectral theorem.

Unlike the authors mentioned above, we shall present a simplified approach, derived directly from the real case, and leading to a more natural formula of the associated spectral representations. 
Specifically, following the author's ideas from \cite{Vas2} (see also \cite{Vas3}), we treat them in a more classical manner, also acting in Hilbert spaces with 
quaternionic bimodule structure, but  where the inner products are real or complex valued, rather than quaternionic valued, as used by most of the previous  authors. In fact, the  quaternionic algebra $\Hh$, and the space of square integrable $\Hh$-valued functions as well, may be endowed with a natural real scalar product. 
Moreover, the spectra and the associated measures of quaternionic normal operators may be  defined in the complex plane, rather than in the quaternionic algebra $\Hh$. Examples in this sense are given in 
the last section of this work.  

The main results of this work in the quaternionic case  are Theorem \ref{Qbdfc} for the bounded case, and 
Theorem \ref{Qubdfc} for the unbounded one. It should be stressed that these results are more or less special cases and direct consequences of Theorems \ref{bdfc0} and  \ref{ubdfc}, respectively. We should  mention that, in the quaternionic case, we also have some uniqueness statements concerning the corresponding spectral measures
(see Propositions \ref{Quniqueness} and  \ref{Qubd_uniqueness}), derived  directly from the results proved in the real case.

We end this work with two  examples of  quaternionic normal operators,
defined by the multiplication with the 
independent quaternionic variable.

Concerning the case of real normal operators, the starting point of our research  was the paper by  Oreshina, who made some useful remarks on a first version of this work. Thanks are also due to an anonymous reader, who suggested the author to mention other contributions, as for instance those due to Goodrich,  Agrawal $\&$ Kulkarni, and Wong, whose  approaches are different from those of this paper. 

We finally mention that, unlike in \cite{AgKu} or in other works, we have constantly use the complex conjugation as a topological involution, thinking that this is the most significant  situation. 
It is clear that replacing the complex conjugation by a topological involution (see for instance \cite{AgKu} for a formal definition and a constant use as well) such arguments can be recaptured via some minor changes.

\section{Normal Operators in Real Hilbert Spaces}
 
\subsection{Preliminaries} 

We start by recalling some elementary facts concerning the linear operators in real Hilbert spaces (see for instance \cite{BaZa},\cite{Ore}, \cite{Vas2} etc.).

Let $\H$ be a real Hilbert space with the inner product $\langle *, *\rangle$, and corresponding norm $\Vert*\Vert$, and let $\H_\C$ be the 
the complexification of $\H$, identified with the direct sum $\H+i\H$.  The natural inner product of $\H_\C$ is then given by
$$
\langle x+iy, u+iv\rangle_\C=\langle x, u\rangle+\langle y, v\rangle +i\langle y, u\rangle-i\langle x, v\rangle
$$
for all $x,y,u,v\in\H$. 

Let also $C:\H_\C\mapsto\H_\C$ be the natural  conjugation $x+iy\mapsto x-iy,\,x,y\in\H$, which is an $\R$-linear isomorphism of $\H_\C$, whose square is the identity.  

We denote by $\mathcal{B}(\H)$ the real algebra of all bounded $\R$-linear operators, acting on $\H$.  Similarly, $\mathcal{B}(\H_\C)$ denotes the complex algebra consisting of all bounded $\C$-linear operators,  acting on $\H_\C$.  

 Each operator 
$T\in\mathcal{B}(\H)$ has a natural extension to an operator
$T_\C\in\mathcal{B}(\H_\C)$, given by $T_\C(x+iy)=Tx+iTy,\,x,y\in\H$.
Moreover, the map $\mathcal{B}(\H)\ni T\mapsto T_\C\in\mathcal{B}(\H_\C)$ is  unital, $\R$-linear, injective,  and multiplicative. In particular, $T\in\mathcal{B}(\H)$ is invertible if and only if $T_\C\in\mathcal{B}(\H_\C)$ is invertible.  The operator $T_\C$ will be 
sometimes called the {\it complex extension} of the real operator $T$.

Fixing an operator $S\in\mathcal{B}(\H_\C)$, we define 
the operator $S^\flat\in\mathcal{B}(\H_\C)$ to be equal to $CSC$.  It is easily seen that the map 
$\mathcal{B}(\H_\C)\ni S\mapsto S^\flat\in \mathcal{B}(\H_\C)$ is a unital conjugate-linear automorphism, whose square is the identity on $\mathcal{B}(\H_\C)$. 
Because $\H=\{u\in\H_\C; Cu=u\}$, we have $S^\flat=S$ if and only if $S(\H)\subset\H$.
In this case, $S=(S\vert\H)_\C$, and so  $T_\C^\flat=T_\C$ for all $T\in\B(\H)$. 

In fact, because 
where $(S+S^{\flat})(\H)\subset\H, i(S-S^{\flat})(\H)\subset\H$,  the algebras $\mathcal{B}(\H_\C)$ and $\mathcal{B}(\H)_\C$ are isomorphic and they may be identified.  In other words, the algebra $\mathcal{B}(\H)$ (in fact $\{T_\C;T\in\mathcal{B}(\H)\}$)  may  be regarded as a (real) subalgebra of $\mathcal{B}(\H_\C)$. In particular, if $S=U+iV$, with $U,V\in\mathcal{B}(\H)$, we have
$S^\flat=U-iV$, so the map $S\mapsto S^\flat$ is the conjugation
of the complex algebra $\mathcal{B}(\H_\C)$ induced
by the conjugation $C$ of $\H_\C$. 

For an operator $T\in\mathcal{B}(\H)$ acting on the real Hilbert space 
$\H$, we denote by $T^*$ its (Hilbert space) adjoint, defined, as usually, via the equality 
$\langle T^*x,y\rangle=\langle x,Ty\rangle$ for all $x,y\in\H$. As we have 
$$
\langle C(x+iy),u+iv\rangle_\C=\langle C(u+iv),x+iy\rangle_\C,
$$
for the operator $T\in
\mathcal{B}(\H)$ it is easily seen that  $(T_\C)^*=(T^*)_\C$, and this operator will be simply denoted by $T^*_\C$. Moreover, $(T_\C^*)^\flat=T_\C^*$.

For every operator $S\in\mathcal{B}(\H_\C)$, we denote,
as usually, by $\sigma(S)$  its spectrum.  Using a classical idea
going back to Kaplansky \cite{Kap} (see also \cite{Ing}), for a real
operator $T\in\mathcal{B}(\H)$ its (complex) spectrum is defined by the equality 
$$
\sigma_\C(T)=\{u+iv;(u-T)^2+v^2\,\, {\rm is\,\,not\,\,invertible}, u,v\in\R\},
$$ 
where the scalars are identified with the corresponding 
multiples of the identity. Because the operator $(u-T)^2+v^2$ is invertible if and only if the operator  
$(u-T_\C)^2+v^2$ is invertible, which in turn is invertible if and only if $u+iv-T_\C$ is
invertible, it follows that $\sigma_\C(T)=\sigma(T_\C)$, as noticed in \cite{Ore} or in \cite{Vas2}. Note also that the set  $\sigma_\C(T)$ is {\it conjugate symmetric}, that is 
$z\in\sigma_\C(T)$ if and only if $\bar{z}\in\sigma_\C(T)$. 

A spectral theory in real Banach spaces has been developed actually 
in the framework of  linear relations in \cite{BaZa}.

Of course, as in the complex case,  a bounded linear operator $T$ in a real  Hilbert space $\H$ is said to be normal if $TT^*=T^*T$. It is easily seen that $T$ is normal if and only if $T_\C$ is normal.

We shall also work with some $\R$-linear operators, not necessarily bounded, called as usually unbounded operators. Such an operator $T$, which is  defined on a vector subspace  $D(T)\subset\H$, with values in $\H$,  has a complex extension $T_\C$, defined on $D(T_\C)=D(T)_\C$, with values in $\H_\C$. When $T$ is closed and/or 
densely defined, then the complex extension $T_\C$ is also closed and/or 
densely defined. 
 
The family of closed and densely defined linear operators  on $\H$,
or  $\H_\C$, will be denoted by $\mathcal{C}(\H)$, or  $\mathcal{C}(\H_\C)$,
respectively.

The adjoint of an operator $T\in\mathcal{C}(\H)$ is defined as  for the complex operators. Adapting the corresponding definition valid in complex Hilbert spaces (as done in
\cite{Ore}), an operator $T\in\mathcal{C}(\H)$
is said to be {\it normal} if $D(T^*)=D(T)$, $D(TT^*)=D(T^*T)$, and $TT^*=T^*T$. 

As in the bounded case,  the complex extension $T_\C$  is normal if and only if
$T$ is normal. Moreover, we have $(T_\C)^*=(T^*)_\C$,  which will be denoted  by $T_\C^*$. 

The spectrum  $\sigma(T_\C)$ of the operator $T_\C\in\mathcal{C}(\H_\C)$ is defined, following \cite{Rud}, as the complement of the set of those points 
$\zeta\in\C$ such that the operator $\zeta-T_\C$ has an everywhere defined bounded inverse. As noticed in \cite{Ore},  the set  
$\sigma(T_\C)$ is still conjugate symmetric. Because we have $T\in\mathcal{C}(\H)$ invertible if and only if $T_\C\in\mathcal{C}(\H_\C)$ is invertible,
and we have $(T_\C)^{-1}=(T^{-1})_\C$, simply denoted by  $T_\C^{-1}$,
we may define  $\sigma_\C(T)=\sigma(T_\C)$, which is the  (complex) spectrum of the real operator $T\in\mathcal{C}(\H)$. Of course, this definition is compatible with the corresponding one in the bounded case. 
  
\subsection{Real and Complex Spectral Measures}

In this subsection we recall some definitions, present some simple consequences of them, and fix some notation. First of all, if $A$
is an arbitrary subset of the complex plane, we put $A^c=\{\bar{z};z\in A\}$. Whe $A=A^c$, we say that $A$ is {\it conjugate symmetric}, a concept already used in the previous subsection.  


\begin{Rem}\label{varprop}\rm (1) In what follows, we shall often use the concept of {\it spectral measure} (see \cite{DuSc}, Part III), which is a particular case of the concept of {\it  resolution of the identity} (see \cite{Rud}), but we shall not distinct the former from the latter, designating them as {\it spectral measures}. 

For our purpose, it is sufficient to consider only spectral measures supported by closed sets in the complex plane. More precisely, a spectral measure in this text is a map
defined on a $\sigma$-algebra, say  ${\bf\Sigma}(\mathfrak{S})$, consisting of Borel subsets of $\mathfrak{S}$, where $\mathfrak{S}$ is a closed subset of the complex plane. Specifically, fixing a real or a complex Hilbert space $\mathcal{K}$,  with the inner product
$\langle *, *\rangle_\mathcal{K}$,  a spectral measure on ${\bf\Sigma}(\mathfrak{S})$  is a set function $F$ on  ${\bf\Sigma}(\mathfrak{S})$ , whose values are self-adjoint projections on   $\mathcal{K}$, with $F(\emptyset)=0$,\, $F(\mathfrak{S})$  the identity of $\mathcal{K},\,F(A\cap B)=F(A)F(B)$ for all $A,B\in{\bf\Sigma}(\mathfrak{S})$, and  such that  the maps $\nu_{x,y}$ on ${\bf\Sigma}(\mathfrak{S})$, given by  $\nu_{x,y}(A)=\langle F(A)x,y\rangle_\mathcal{K}$ for all $x,y\in \mathcal{K}$, are countably additive measures.  The scalar measures $\nu_{x,y}$, called the {\it associated scalar measures} to $F$, are complex  (resp. real) if the space  $\mathcal{K}$ is complex (resp. real; for the real case see also \cite{Ore}). 
In particular, if $\{A_k\}_{k=1}^\infty$ is a countable family of mutually disjoint sets from ${\bf\Sigma}(\mathfrak{S})$, we have $F(\cup_{k=1}^\infty A_k)x=\sum_{k=1}^\infty F(A_k)x$ for every $x\in \mathcal{K}$.
\medskip

(2) If $\H$ is a real Hilbert space, considering its complexification $\H_\C$, and fixing 
a spectral measure  $E: Bor(\mathfrak{S})\mapsto\mathcal{B}(\H_\C)$,   
we may consider its restriction $E_\R(A)=E(A)\vert\H$ for all $A\in Bor_E(\mathfrak{S})$, where $ Bor_E(\mathfrak{S})=\{A\in Bor(\mathfrak{S}); E(A)=E(A)^\flat\}$. The map
 $E_\R: Bor_E(\mathfrak{S})\mapsto\mathcal{B}(\H)$ is a real spectral measure, because ${\rm Bor}_E(\mathfrak{S})$ is a $\sigma$-algebra. Indeed, we have $\mathfrak{S}\setminus A\in {\rm Bor}_E(\mathfrak{S})$ whenever $A\in Bor_E(\mathfrak{S})$, and $E_\R(\cup_{k=1}^\infty A_k)=
E_\R(\cup_{k=1}^\infty A_k)^\flat$  if $\{A_k\}_{k=1}^\infty$ is a sequence of mutually disjoint sets from $Bor_E(\mathfrak{S})$, via the additivity and continuity of the map $\mathcal{B}(\H_\C)\ni S\mapsto S^\flat\in \mathcal{B}(\H_\C)$. The general case of an arbitrary sequence $\{B_k\}_{k=1}^\infty$   from $Bor_E(\mathfrak{S})$ can be reduced to the previous one, replacing it 
by a sequence of mutually disjoint sets, having the same union.

Note also that if $\{\nu_{\xi,\eta};\xi,\eta\in\H\}$ is the family of
associated scalar measures of $E$, then $\{\nu_{x,y};x,y\in\H\}$ is a 
family of associated scalar measure of $E_\R$ actually defined on 
the whole $\sigma$-algebra  $Bor(\mathfrak{S})$.
\medskip

(3) We denote by $\mathfrak{B(S)}$ the complex algebra of all complex-valued Borel
functions, defined on $\mathfrak{S}$. When $\mathfrak{S}$ is conjugate symmetric,
fixing a complex Hilbert space  $\mathcal{K}$,  and  a spectral measure 
$E:Bor(\mathfrak{S})\mapsto\mathcal{B}(\mathcal{K})$,
we denote by $\mathfrak{B_{\it s,E}(S)}$ the real subalgebra of those  functions  $f\in\mathfrak{B(S)}$ with the property $f(\bar{z})=\overline{f(z)}$ for all $z\in \mathfrak{S}$, except for a set of null $E$-measure, which will be sometimes designated as  
$E$-{\it stem functions}.
\medskip

(4) As in \cite{Rud}, Section 12.20, fixing again a complex Hilbert space  $\mathcal{K}$,  and  a spectral measure  $E: Bor(\mathfrak{S})\mapsto\mathcal{B}(\mathcal{K})$,  we consider the algebra ${\mathcal L}^\infty(\mathfrak{S},E)$  consisting of  $E$-essentially bounded measurable functions  $f\in\mathfrak{B(S)}$, that is,
$$
\Vert f\Vert_\infty=\inf\{r>0;E(\{z\in\mathfrak{S},\vert f(z)\vert \ge r\})=0\}<\infty.
$$
Setting $\mathcal{N}=\{f\in{\mathcal L}^\infty(\mathfrak{S},E);\Vert f\Vert_\infty=0\} $, which is un ideal in ${\mathcal L}^\infty(\mathfrak{S},E)$,
we define the quotient $L^\infty(\mathfrak{S},E)=
{\mathcal L}^\infty(\mathfrak{S},E)/\mathcal{N}$, which is a Banach algebra.
In fact, $L^\infty(\mathfrak{S},E)$ is a $C^*$-algebra, with the involution induced by 
the complex conjugation $f\mapsto\bar{f}$, where $\bar{f}(z)=\overline{f(z)},\,
z\in\mathfrak{S}$.
 
Practically, we identify a function $f\in{\mathcal L}^\infty(\mathfrak{S},E)$ with its 
equivalence class $f+\mathcal{N}\in L^\infty(\mathfrak{S},E)$.

In particular, given two sets $A,B\in Bor(\mathfrak{S})$, by abuse of 
notation, we sometimes  write $A=B$
if their characteristic functions $\chi_A,\chi_B$ are in the same 
equivalence class, that is $\chi_A=\chi_B$ in $L^\infty(\mathfrak{S},E)$. This is also equivalent to the equalities
$E(A\setminus B)=0=E(B\setminus A)$. 
 
(5) Let again $\H$ to be a real Hilbert space, take its complexification $\H_\C$, and fix a spectral measure  $E:{\rm Bor}(\mathfrak{S})\mapsto\mathcal{B}(\H_\C)$.

 Assuming $\mathfrak{S}$ conjugate symmetric, we denote by
 $L^\infty_s(\mathfrak{S},E)$ the (real) subalgebra of  
$L^\infty(\mathfrak{S},E)$  consisting of $E$-stem functions.

Note that the finite sums of the form $\Sigma_{j\in J}(r_j\chi_{A_j}+is_j\theta_{A_j})$, with
$r_j,s_j\in\R$ for all $j\in J,\,(A_j)_{j\in J}\subset Bor_E(\mathfrak{S})$  a partition of $\mathfrak{S}$,  where $\chi_A$ is the characteristic function of the set $A$, and $\theta_A(z)=1,-1,=0$ if $\Im(z)>0$
$\Im(z)<0$, and either  $\Im(z)=0$ or  $z\in\mathfrak{S}\setminus A$, respectively, are elements of $L^\infty_s(\mathfrak{S},E)$. The functions of this type will be designated in this text as 
{\it elementary functions} in  
$ L^\infty_s(\mathfrak{S})$, which denotes an abbreviation of  $L^\infty_s(\mathfrak{S},E)$.

Let  $F^\infty_{s}(\mathfrak{S})$ be the real subspace of $L^\infty_s(\mathfrak{S})$,
generated by  the elementary functions. We can easily see that the subspece $F^\infty_{s}(\mathfrak{S})$ is dense in the algebra  $L^\infty_s(\mathfrak{S})$. Fixing a  function $f\in L^\infty_s(\mathfrak{S})$, we 
can represent it as $f(z)=f_1(z)+if_2(z)$, with $f_1,f_2$ real valued, and with 
$f_1(\bar{z})=f_1(z)$, and $f_2(\bar{z})=-f_2(z)$.

Taking the restriction $g_1(z)=f_1(z)$ to the set  $\{z\in \mathfrak{S};\Im(z)\ge0\}$, we can find a simple function $h$ arbitrarily 
close to the function  $g_1$. Extending $h$ to the whole set  $\mathfrak{S}$ by putting $h(\bar{z})=h(z)$, we obtain an elementary function, say $h_1$, which is arbitrarily close to $f_1$. Similarly, we can construct an elementary function $h_2$, arbitrarily close
to $f_2$, showing eventually that the family of elementary functions is dense in  $L^\infty_s(\mathfrak{S})$.

\end{Rem}


\begin{Rem}\label{dir_sum0}\rm \rm Let $\mathfrak{S}\subset\C$ be a conjugate symmetric set, let $\tau:\mathfrak{S}\mapsto\mathfrak{S}$
be given by $\tau(z)=\bar{z},\,z\in\mathfrak{S}$,
and let $\mathfrak{F}(\mathfrak{S})$ be a vector space of complex-valued functions on $\mathfrak{S}$, with the property
$f\in \mathfrak{F}(\mathfrak{S})$ if and only if $f\circ\tau\in \mathfrak{F}(\mathfrak{S})$. Also set $\bar{f}(z)=\overline{f(z)}$ for all $z\in \mathfrak{S}$. Then we have a direct sum decomposition
$\mathfrak{F}(\mathfrak{S})=\mathfrak{F}_s(\mathfrak{S})+\mathfrak{F}_a(\mathfrak{S})$, where

$$
\mathfrak{F}_s(\mathfrak{S})=\{f\in \mathfrak{F}(\mathfrak{S}); f=\bar{f}\circ\tau\} 
$$
and 
$$
\mathfrak{F}_a(\mathfrak{S})=\{f\in \mathfrak{F}(\mathfrak{S}); f=-\bar{f}\circ\tau\}
$$

It is obvious that $\mathfrak{F}_s(\mathfrak{S})\cap \mathfrak{F}_a(\mathfrak{S})=\{0\}$. Note also that that $\mathfrak{F}(\mathfrak{S})=\mathfrak{F}_s(\mathfrak{S})+\mathfrak{F}_a(\mathfrak{S})$. Indeed, fixing an 
arbitrary function $f$ from $\mathfrak{F}(\mathfrak{S})$, we can write   $f=g+h$, where $g=2^{-1 }(f+ \bar{f}\circ\tau)$ and 
$h=2^{-1 }(f- \bar{f}\circ\tau)$, with 
$g\in \mathfrak{F}_s(\mathfrak{S})$ and 
$h=\mathfrak{F}_a(\mathfrak{S})$.
\end{Rem}
  
  
\begin{Rem}\label{dir_sum}\rm Let $\mathfrak{S}\subset\C$ be a compact conjugate 
symmetric set, and let $C(\mathfrak{S})$ be the algebra of all 
complex valued continuous functions on $\mathfrak{S}$. As in Remark \ref{dir_sum0},
considering  the subspaces 
$$
C_s(\mathfrak{S})=\{f\in C(\mathfrak{S}); f=\bar{f}\circ\tau\} 
$$
and 
$$
C_a(\mathfrak{S})=\{f\in C(\mathfrak{S}); f=-\bar{f}\circ\tau\}, 
$$
we have $C(\mathfrak{S})=C_s(\mathfrak{S})+C_a(\mathfrak{S})$, and the 
sum is direct. 

A classical result due to Riesz asserts the dual space $C(\mathfrak{S})^*$ of $C(\mathfrak{S})$ can be identified with the space $M(\mathfrak{S})$ consisting of all Borel measures on $\mathfrak{S}$. 
The dual of the space $C_s(\mathfrak{S})^*$ can also be described in terms of measures (see \cite{Grz}), and it can be identified with the space $M_s(\mathfrak{S})$ consisting of all measures $\nu$ on $Bor(\mathfrak{S})$ with the property $\nu=\bar{\nu}\circ\tau$.
\end{Rem}


\begin{Rem}\label{dir_sum1}\rm Let, as above,  
$\mathfrak{S}\subset\C$ be a compact conjugate symmetric open set. Fixing a complex Hilbert space  $\mathcal{K}$,  and  a spectral measure  $E: Bor(\mathfrak{S})\mapsto\mathcal{B}(\mathcal{K})$,  we have introduced the algebra 
$ L^\infty(\mathfrak{S},E)$. Again as in As in Remark \ref{dir_sum0},
we have that $L^\infty(\mathfrak{S},E)$ is the direct sum of the spaces $L^\infty_s(\mathfrak{S},E)+L^\infty_a(\mathfrak{S},E)$, where we recall
that
$$
L^\infty_s(\mathfrak{S},E)=\{f\in L^\infty(\mathfrak{S},E); f(z)=\overline{f(\bar{z})},
z\in\mathfrak{S}\},
$$
and we put
$$
L^\infty_a(\mathfrak{S},E)=\{f\in L^\infty(\mathfrak{S},E); f(z)=-\overline{f(\bar{z})}, z\in\mathfrak{S}\}. 
$$
Therefore, there  exists a real linear projection, say $P_\infty$,
form $L^\infty(\mathfrak{S},E)$ onto $L^\infty_s(\mathfrak{S},E)$
given by $P_\infty f=2^{-1}(f+\bar{f}\circ\tau)$.

Note that $\chi_A\in L^\infty_s(\mathfrak{S},E)$ if and only if 
$\chi_A=\chi_{A^c}$, where 
$\chi_A$ is the characteristic function of the set $A$. Moreover, if $A\in Bor(\mathfrak{S})$ is arbitrary, then $P_\infty\chi_A=
2^{-1}(\chi_A+\chi_{A^c})$.
\end{Rem}


\subsection{Bounded Normal Operators in Real Hilbert Spaces}

Following the corresponding part of \cite{Rud} (see  the subsections 12.17-12.26), we recall some general results concerning the bounded normal operators, necessary for our further development.


\begin{Rem}\label{cst}\rm Let $\mathcal{K}$ be a complex Hilbert space, and let $S\in\mathcal{B(K)}$ be a 
normal operator. If $\mathfrak{S}$ is the spectrum of $S$, there exists a unique 
spectral measure $E: Bor(\mathfrak{S})\mapsto\mathcal{B(K)}$ inducing a
unital $C^*$-algebra morphism
$$
L^\infty(\mathfrak{S},E)\ni f\mapsto f(S):= \int_\mathfrak{S}f(z)dE(z)\in \mathcal{B(K)},
$$
meaning that 
$$
\langle f(S)x,y\rangle_\mathcal{K}=\int_{\mathfrak{S}}f(z) d\mu_{x,y}(z),\,\,
x,y\in\mathcal{K},
$$
where $\{\mu_{x,y};x,y\in\mathcal{K}\}$ are the scalar measures associated to $E$.

We are especially interested to obtain a spectral representation of 
a normal operator $T\in\mathcal{B(H)}$, in the case of a real Hilbert space $\H$, using known results concerning its  
normal extension $T_\C\in\mathcal{B(H}_\C)$, whose  spectral measure 
$E_\C: Bor(\mathfrak{S})\mapsto\mathcal{B(H}_\C)$ is supposed to be given, where $\mathfrak{S}=\sigma(T_\C)$. As the real Hilbert space  $\H$ is not necessarily invariant under the spectral measure $E_\C(*)$,  we shall charactrize those Borel subsets $A\subset\mathfrak{S}$ such that $E_\C(A)^\flat=E_\C(A)$. 
\end{Rem}

We recall that for an arbitrary subset $A\subset\C$, we put $A^c:=\{\bar{z};z\in A\}$. 


\begin{Lem}\label{inv_cond} With $T,\,T_\C$ and $E_\C$ as above,  we have the equality $E_\C(A)^\flat=E_\C(A)$ for some $A\in{\rm Bor}(\mathfrak{S})$ if and only if $\chi_A=\chi_{A^c}$ 
in $L^\infty(\mathfrak{S},E_\C)$. 
\end{Lem}

{\it Proof}\,  We first note that the map $E_\C^\flat:{\rm Bor}(\mathfrak{S})\mapsto\mathcal{B}(\H_\C)$, with $E_\C^\flat(A)=E_\C(A)^\flat$, is also a
spectral measure, which is easily checked. Let us find the operator $S=\int_\mathfrak{S} z dE_\C^\flat(z)$. Note that for all 
$\xi,\eta\in\H_\C$, and $A\subset\mathfrak{S}$ a Borel subset, we have
$$
\langle E_\C^\flat(A)\xi,\eta\rangle_\C=\overline{\langle  E_\C(A)C\xi,C\eta\rangle_\C},
$$
and so $\mu^\flat_{\xi,\eta}=\bar{\mu}_{C\xi,C\eta}$, where $\mu^\flat_{\xi,\eta}$
is the scalar measure  $\langle E_\C^\flat(*)\xi,\eta\rangle_\C$. Therefore,
$$
\langle S\xi,\eta\rangle_\C=\int_\mathfrak{S} z d\mu^\flat_{\xi,\eta}=\int_\mathfrak{S} z d\bar{\mu}_{C\xi,C\eta}=
\overline{\int_\mathfrak{S} \bar{z} d\mu_{C\xi,C\eta}}
$$
$$
=\overline{\langle T_\C^*C\xi,C\eta\rangle_\C}=\langle  CT_\C^*C\xi,\eta\rangle_\C=
\langle  T_\C^*\xi,\eta\rangle_\C,
$$
showing that  $E_\C^\flat$ is precisely the spectral measure of the normal operator
$T^*_\C$. 

The spectral measure of the adjoint  $T^*_\C$ of the normal operator $T_\C$
can be also obtain by regarding   $T^*_\C$ as a function of $T_\C$, via the map
$\C\ni z\mapsto\bar{z}\in\C$. Applying a particular case of the change of measure
principle (see \cite{Rud}, Theorem 13.28), we obtain the equality  $E_\C^\flat(A)= E_\C(A^c)$ for all $A\in{\rm Bor}(\mathfrak{S})$. 

In particular, $E_\C(A)=0$ if and only if $E_\C(A^c)=0$.

If $\chi_A=\chi_{A^c}$, we must have  have $E_\C(A)^\flat=E_\C(A)$. Conversely,
assume that $E_\C(A)^\flat=E_\C(A)$. Writing $A=(A\setminus A^c)\cup(A\cap A^c)$ and $A^c=(A^c\setminus A)\cup(A\cap A^c)$, we deduce that $E_\C(A^c\setminus A)=E_\C(A\setminus A^c)$. On the other 
hand, the product $E_\C(A\setminus A^c)E_\C(A^c\setminus A)=E_\C(\emptyset)=0$.
Therefore, $E_\C(A\setminus A^c)=E_\C(A^c\setminus A)=0$ implying
$\chi_A=\chi_{A^c}$ in  $L^\infty(\mathfrak{S}, E_\C)$, which completes the proof. 

\noindent{\bf QED}


\begin{Rem}\label{inv_cond1}\rm (1) Lemma \ref{inv_cond} establishes a partial  invariance of the real Hilbert space $\H$, with respect to the spectral measure $E_\C$ of the normal operator $T_\C\in\mathcal{B}(\H_\C)$, obtained via a normal operator $T$, acting on $\H$. In fact, via Remark \ref{cst}
we have that the operator $\chi_A(T_\C)=E_\C(A)$ leaves invariant the space $\H$ 
if and only if $\chi_A=\chi_{A^c}$ in $L^\infty(\mathfrak{S},E_\C)$, written as $A=AZ^c$, that is,   when $E_\C(A^c\setminus A)=E_\C(A\setminus A^c)=0$ (as in Remark \ref{varprop} (4)).  

Further, setting $E_\R(A)=E_\C(A)\vert \H=\chi_A(T_\C)\vert\H$, we get a family of self-adjoint projections on $\H$. Moreover, with our convention, 
$$\{A\in{Bor}(\mathfrak{S});E_\C(A^c)=E_\C(A)\}=\{A\in{Bor}(\mathfrak{S});A^c=A\},
$$
is a $\sigma$-algebra denoted by $Bor_{E_\C}(\mathfrak{S})$ (see Remark \ref{varprop}(2)). Moreover, the map
$Bor_{E_\C}(\mathfrak{S})\ni A\mapsto E_\R(A)\in \mathcal{B}(\H)$ 
is a real spectral measure. Some properties of the real spectral measure $ E_\R$ can be deduced from the properties of the spectral measure $E_\C$. 
\medskip

(2) Simple functions of the form $\sum_{j\in J} r_j A_j$ 
with $r_j\in\R$ for all $j\in J$, and with $(A_j)_{j\in J}\subset Bor_{E_\C}(\mathfrak{S})$  a partition of $\mathfrak{S}$ may be integrated with
respect to $E_\R$, leading to the operator 
$\sum_{j\in J} r_j E_\R(A_j)$. Nevertheless, the spectral measure $E_\R$
does not seem to be 'strong` enough to be used to integrate  other elementery functions. For instance,  the function 
$h=i(\chi_A-\chi_{A^c})$, where $A\in Bor(\mathfrak{S})$ is arbitrary, cannot be directly integrated with respect to the measure $E_\R$, when $A\neq A^c$.  Nevertheless, the space $\H$ is invariant under the operator  
$h(T_\C)=i(E_\C(A)-E_\C(A^c))$, which is the integral of $h$ 
with respect to $E_\C$ because $h(T_\C)^\flat=h(T_\C)$. Hence, we may define $h(T)=h(T_\C)\vert\H$. 

A  related discussion can be found in Example \ref{con_exa}.
\end{Rem}


\begin{Rem}\label{rGt}\rm Let $T\in\mathcal{B}(\H)$ be a real normal operator, and let $T_\C\in \mathcal{B}(\H_\C)$ be its complex extension.
We may consider in $\B(\H_\C)$ the unital complex Banach algebra 
$\A_\C$, generated by the operators $T_\C, T_\C^*$, which is, in fact, a commutative complex $C^*$-algebra.
  
Let us denote by $\A_\R$ the real subalgebra of $\A_\C$, 
consisting  of limits of sequences of real linear combinations of operators $T_\C^m T_\C^{*n}$, with $m,n$ nonnegative integers, which is  a real $C^*$-algebra. The elements of $\A_\R$ are operators which leave invariant the space $\H$.   

According to Gelfand's theory, the algebra $\A_\C$ may be identified with the algebra $C(\mathfrak{S})$ of all complex valued continuous functions on $\mathfrak{S}$. For the real algebra $\A_\R$ there is also a Gelfand type theory, due to 
Arens and Kaplansky (see \cite{Good}), allowing the identification of  this algebra with the real algebra $C_s(\mathfrak{S})$  of all functions $f\in C(\mathfrak{S})$ with the
property $f(z)=\overline{f(\bar{z})}$ for every $z\in\mathfrak{S}$
(see Remark $\ref{dir_sum}$). Also note that the polynomials $
p(z,\bar{z})$ with real coefficients are uniformly dense in $C_s(\mathfrak{S})$ via  a real analogue of the Stone-Weierstrass Theorem (see \cite{KulLim}). For this reason, the algebra $C_s(\mathfrak{S})$ plays an important role in the study of real normal
operators.  
\end{Rem}

The next statement is a representation theorem, in fact a functional calculus for real normal operators with  functions from the real
algebra $L^\infty_s(\mathfrak{S},E_\C)$ (see Remark \ref{dir_sum1}), 
obtained by using some results from \cite{Rud} (especially Theorem 12.21).

In the following, we designate by $[A]$ the closure of the arbitrary subset $A\subset\C$.

\bigskip

\begin{Thm}\label{bdfc0} Let $T\in\mathcal{B}(\H)$ be a real normal operator, let $T_\C\in \mathcal{B}(\H_\C)$ be its complex extension, and let $E_\C:Bor(\mathfrak{S})
\mapsto \mathcal{B}(\H_\C)$ be the spectral measure of $T_\C$, where
$ \mathfrak{S}=\sigma(T_\C)$.

We consider the map $\Phi:L^\infty_s (\mathfrak{S},E_\C)
\mapsto\mathcal{B}(\H)$,
given by the formula
$$
\langle\Phi(f)x, y \rangle_\H=\int_\mathfrak{S} fd\mu_{x,y},\,\, f\in L^\infty_s (\mathfrak{S},E_\C),\, x,y\in\H,
$$
where $\{\mu_{x,y};x,y\in\H\}$ is the corresponding subset of the family of   measures 
$\{\mu_{\xi,\eta}; \xi,\eta\in\H_\C\}$ associated with the spectral measure $E_\C$. 

The map $\Phi$ 
is a unital real  algebra isometric morphism from $ L^\infty_s (\mathfrak{S},E_\C)$ into $\mathcal{B}(\H)$.  Moreover, we have the following properties.

(i) $\Phi(f)$ is a real normal operator and $\Phi(\bar{f})=\Phi(f)^*=
\Phi(f\circ\tau)$ for all 
$f\in L^\infty_s (\mathfrak{S},E_\C)$, where $\tau(z)=\bar{z}$ for all
$z\in\mathfrak{S}$.

(ii) For every polynomial $p(z,\bar{z})$ with real coeficients one has 
$\Phi(p)=p(T,T^*)$.

(iii) $\Vert \Phi(f)x\Vert^2_\H=\int_\mathfrak{S}\vert f\vert^2d\mu_{x,x}$ for all 
$f\in L^\infty_s (\mathfrak{S},E_\C)$ and $x\in\H$. 

(iv) We have the inclusion 
$\sigma(\Phi(f))\subseteq [f(\mathfrak{S})]$
for all $f\in L^\infty_s (\mathfrak{S},E_\C).$

(v)  If  $S\in\mathcal{B}(\H)$ and $ST=TS$, then $S\Phi(f)=\Phi(f)S$ 
 for all $f\in L^\infty_s (\mathfrak{S},E_\C)$.
\end{Thm}

{\it Proof}\,  As in Remark \ref{cst}, we have a $C^*$-algebra morphism from $L^\infty(\mathfrak{S},E_\C)$ into $\B(\H_\C)$,  explicitly written as
$$
\langle f(T_\C)\xi,\eta\rangle_{\H_\C}=\int_{\mathfrak{S}}f(z) 
d\mu_{\xi,\eta}(z),\,\,\xi,\eta\in\H_\C,\,f\in L^\infty(\mathfrak{S},E_\C).
$$
Let us show the space $\H$ is invariant under the operator $f(T_\C)$ for all 
$f\in L^\infty_s (\mathfrak{S},E_\C)$. For the elements of the algebra
$\A_\R$, this property is already mentioned in Remark \ref{rGt}. 
Nevertheless, we need a stronger result concerning the real algebra 
$ L^\infty_s(\mathfrak{S},E_\C)$. 

Fixing  $A\in\rm{Bor}_{E_\C}(\mathfrak{S})$, we set $A^+=\{z\in A;\Im{z}>0\},\,A^-=\{z\in A;\Im{z}<0\},\,A^0=\{z\in A;\Im{z}=0\}$. Clearly, ${A^+}^c=A^-$.
Moreover, with $\theta_A$ as in Remark \ref{varprop}(5),
$$
\int_\mathfrak{S} \theta_A dE_\C=E_\C(A^+)- E_\C(A^-),
$$
and $C[i(E_\C(A^+)- E_\C(A^-))]C=i(E_\C(A^+)- E_\C(A^-))$, showing that the space $\H$ is invariant under the operator $i\theta_A(T_\C)=i(E_\C(A^+)- E_\C(A^-))$ (see also Remark \ref{inv_cond1}(2)).

Choosing $f=\Sigma_{j\in J}(r_j\chi_{A_j}+is_j\theta_{A_j})$, that is,  an elementary function, and for $x,y\in\H$, we have
$$
\int_\mathfrak{S} \Sigma_{j\in J}(r_j\chi_{A_j}+is_j\theta_{A_j})d\mu_{x,y}=
$$
$$
\langle \Sigma_{j\in J} (r_jE_\C(A_j)+is_j\theta_{A_j}(T_\C))x,y\rangle_\H,
$$
hence $\Phi(f)=  \Sigma_{j\in J}(r_jE_\C(A_j)+is_j\theta_{A_j}(T_\C))\vert\H=f(T_\C)\vert\H$, using also Lemma \ref{inv_cond}.

Let  $L^\infty_{s0}(\mathfrak{S})$ be the subalgebra of $L^\infty_s(\mathfrak{S})=L^\infty_s(\mathfrak{S},E_\C)$,
generated by  the elementary functions. Because the map $L^\infty(\mathfrak{S})\ni f\mapsto f(T_\C)\in\mathcal{B}(\H_\C)$ is an algebra morphism, it follows, in particular, that for every $f\in L^\infty_{s0}(\mathfrak{S})$, the space $\H$ is invariant under 
$ f(T_\C)$, and so $\Phi(f)= f(T_\C)\vert\H_\C$ is $\H$-valued.
Hence,  using the density of $ L^\infty_{s0}(\mathfrak{S})\supset  F^\infty_{s}(\mathfrak{S})$ in $L^\infty_s(\mathfrak{S})$ (see Remark \ref{varprop}(5)), a direct  continuity argument shows that   the space $\H$ is invariant under $ f(T_\C)$ for all $f\in L^\infty_s(\mathfrak{S})$.

In fact,  the map 
$\Phi(f):L^\infty_s (\mathfrak{S})\mapsto\mathcal{B}(\H)$ is the restriction 
of the map $ L^\infty(\mathfrak{S})\ni f\mapsto f(T_\C)\in\mathcal{B}(\H_\C)$, which is
the functional calculus of $T_\C$ with bounded masurable functions.
Therefore, it inherits most of the properties of the latter as stated 
in \cite{Rud} (see especially sections 12.20-12.24).

In particular, it is an isometric morphism of real 
algebras, whose image consists of normal operators. When 
$f\in L^\infty_s (\mathfrak{S},E_\C)$,
we must have $(i)$, the second equality following from the property
$\bar{f}=f\circ\tau$. The properties $(ii),(iii)$ are also direct consequences from the properties of the map $ f\mapsto f(T_\C)$.

The inclusion $(iv)$ is obtained in the following way:
$$
\sigma(\Phi(f))=\sigma(\Phi(f)_\C)=\sigma(f(T_\C))\subseteq [f(\sigma(T_\C))]=
[f(\mathfrak{S})]
$$
for all $f\in L^\infty_s (\mathfrak{S},E_\C)$, which is a particular case of  Corollary X.2.9(III) from 
\cite{DuSc}, Part II. (See also our Remark \ref{re_spe_me} for a more general assertion.)

If $S\in\mathcal{B}(\H)$ and $ST=TS$, we must also have $S_\C T_\C=T_\C S_\C$.
This implies  $S_\C f(T_\C)= f(T_\C)S_\C$ for all $f\in L^\infty_s (\mathfrak{S},E_\C)$, via a classical theorem due tu Fuglede \cite{Fug}.
Passing to restrictions to $\H$, we obtain  $S\Phi(f)=\Phi(f)S$ for all 
 $f\in L^\infty_s (\mathfrak{S},E_\C)$ showing that $(v)$
also holds.

\noindent{\bf QED}


\begin{Rem}\label{re_spe_me}\rm With the notation of Theorem \ref{bdfc0}, we can associate  every real normal operator $T\in\B(\H)$ with a spectral measure $E_\R$, which, according to Lemma \ref{inv_cond}, should be defined on the $\sigma$-algabra 
$Bor_{E_\C}(\mathfrak{S})$ (see Remark 6). It is given by
$$
\langle E_\R(A)x,y\rangle_\H=\int_\mathfrak{S}\chi_A(z)d\mu_{x,y}(z),\,\,
A\in Bor_{E_\C}(\mathfrak{S}),\,\,x,y\in\H,
$$
that is, $E_\R(A)=E_\C(A)\vert\H$ for all $A\in Bor_{E_\C}(\mathfrak{S})$. We also have
 $$
 \langle Tx,y\rangle_\H=\int_\mathfrak{S} z d\mu_{x,y}(z),\,\,x,y\in\H.
 $$
This is a spectral theorem using a real spectral measure 
described in Remark \ref{varprop}(2) (see also Definition 2.1 and Theorem 2.7 from \cite{AgKu}). 

Note that we have 

$$\sigma(\Phi(f))=\bigcap_{E_\R(A)=I_{\H}}[f(A)],$$
where $I_\H$ is the identity on $\H$.

To get this equality, we use the following formula
$$\sigma(f(T_\C))=\bigcap_{E_\C(A)=I_{\H_\C}}[f(A)],$$
proved in Corollary X.2.9(III) from \cite{DuSc}, Part II, 
applied to $f\in L^\infty_s (\mathfrak{S},E_\C)$, where $I_{\H_\C}$ is 
the identity on $\H_\C$.  Note that the equality $E_\C(A)=I_{\H_\C}$ 
implies $E_\C(A)=E_\C(A^c)$ because $E_\C(A^c)=E_\C(A)^\flat=
I_{\H_\C}^\flat=I_{\H_\C}$, via Lemma \ref{inv_cond}. Thus
$E_\C(A)=I_{\H_\C}$ if and only if $E_\R(A)=I_{\H}$.

The inclusion $(iv)$ from Theorem \ref{bdfc0} follows directly from 
this equality.

Note also that 
 $\sigma(T\vert E_\R(A)\H)=\sigma(T_\C\vert E_\C(A)\H_\C)\subseteq[A]$
 for all $A\in Bor_{E_\C}(\mathfrak{S})$, which is precisely
 Corollary X.2.6 from \cite{DuSc}, Part II. 
 
 We may and shall sometimes use the notation $f(T)$ to designate the operator $\Phi(f)$  for every $f\in L^\infty_s (\mathfrak{S},E_\C)$.

We finally remark that a real normal operator $T$ may be also viewed as a real version of the concept of {\it spectral operator} in the sense of Definition XV.2.5 from \cite{DuSc}, Part III, because it can be associated with the real spectral measure $E_\R$,  satisfying the 
properties $TE_R(A)=E_R(A)T$ and $\sigma(T\vert E_R(A)\H)\subset[A]$
for all $A\in Bor_{E_\C}(\mathfrak{S})$, as required by Definition 
XV.2.2 from \cite{DuSc}, Part III. 
\end{Rem}

The following result  insures the uniqueness of the measures 
$\{\mu_{x,y};x,y\in\H\}$ from our Theorem \ref{bdfc0}.


\begin{Cor}\label{bdfc1} Let $T\in\mathcal{B}(\H)$ be a real normal operator, let $T_\C\in \mathcal{B}(\H_\C)$ be its complex extension,  let $E_\C:Bor(\mathfrak{S})
\mapsto \mathcal{B}(\H_\C)$ be the spectral measure of $T_\C$, where
$ \mathfrak{S}=\sigma(T_\C)$, and let $\{\mu_{x,y};x,y\in\H\}$ be the corresponding subset of the family of   measures 
$\{\mu_{\xi,\eta}; \xi,\eta\in\H_\C\}$ associated with the spectral measure $E_\C$. Then there exists a unique family of measures 
$\{\nu_{x,y};x,y\in\H\}$ such that 
$$
\int_{\mathfrak{S}}fd\nu_{x,y}=\int_{\mathfrak{S}}fd\mu_{x,y}=\langle f(T)x,y\rangle_\H, \,f\in C_s(\mathfrak{S}), 
$$ 
with $\nu_{x,y}=\bar{\nu}_{x,y}\circ\tau=\mu_{x,y}$ for all $x,y\in\H$. 
\end{Cor}

{\it Proof}\, Because the real Hilbert space $\H$ is invariant under
$f(T_\C)$ for all $f\in C_s(\mathfrak{S})$, it follows that the linear
map $C_s(\mathfrak{S})\ni f\to\langle f(T_\C)x,y\rangle_\H=
\langle f(T)x,y\rangle_\H$ is a real linear functional on $C_s(\mathfrak{S})$. 

According to the real analogue of the Riesz Representation Theorem due to Grzesiak (see \cite{Grz}), there exists a complex measure 
$\nu_{x,y}$ in the dual of the space $C_s(\mathfrak{S})$   such that 
$$
\int_{\mathfrak{S}}fd\nu_{x,y}=\langle f(T)x,y\rangle_\H =\int_{\mathfrak{S}}fd\mu_{x,y}, \,f\in C_s(\mathfrak{S}),\,x,y\in\H. 
$$ 

Moreover, we must have $\nu_{x,y}=\bar{\nu}_{x,y}\circ\tau$ for all $x,y\in\H$ via  Grzesiak's theorem, also implying the uniqueness
of this family.  

In fact, when applying the equality from above to $f=ig\in C_s(\mathfrak{S})$ with $g\in C_a(\mathfrak{S})$.
we obtain that the equality of the integrals is valid for all $f\in C(\mathfrak{S})$,
insuring the equality of the measures $\nu_{x,y}=\mu_{x,y}$ for all $x,y\in\H$. 

\noindent{\bf QED}

\begin{Exa}\label{con_exa}\rm We recall Example 2.2 from \cite{AgKu}, which will be 
re-examined.  Let $L^2(\mathfrak{S})$ be the space of all complex valued square integrable functions on $\mathfrak{S}$ with respect to the Legesgue measure denoted by $\lambda$, where $\mathfrak{S}$
is the closed unit disk in $\C$. Let also
$$
L^2_s(\mathfrak{S})=\{f\in L^2(\mathfrak{S}); f(z)=\overline{f(\bar{z})},
z\in\mathfrak{S}\},
$$
which is a real Hilbert subspace of $L^2(\mathfrak{S})$ , whose scalar product is the restriction of that of $L^2(\mathfrak{S})$. In othe words, $\langle\phi,\psi\rangle =\int_\mathfrak{S} \phi\bar{\psi}d\lambda$.

In fact, we 
have $L^2(\mathfrak{S})=L^2_s(\mathfrak{S})_\C$. The multiplication 
operator $Tf(z)=zf(z),\,z\in \mathfrak{S}$ is a normal operator, and so
the complex extension $T_\C$ of $T$ is a normal operator on  $L^2(\mathfrak{S})$, whose spectral measure is given by $E_\C(A)\phi=
\chi_A\phi$ for all $A\in Bor(\mathfrak{S})$ and $\phi\in L^2(\mathfrak{S})$. The spectral measure $E_\R$  of $T$ is the restriction of $E_\C$ to $Bor_{E_\C}(\mathfrak{S})=\{A\in Bor(\mathfrak{S};
\lambda(A)=\lambda(A^c)\}$.  Nevertheless,  the attached scalar measures of $E_\R$, when defined only on the $\sigma$-algabra $Bor_{E_\C}(\mathfrak{S})$ cannot integrate some simple functions. Indeed, 
as in Remark \ref{inv_cond1}(2), the function
$h(z)=i(\chi_A(z)-\chi_{A^c}(z)$ from $L^2_s(\mathfrak{S})$, where $\lambda(A)\neq \lambda(A^c)$, cannot be integrated with respect to such a measure. 
\end{Exa}

The next assertion is a uniqueness result for the real spectral measure attached to a real normal operator. Using only (spectral) inclusions, this result is different (and more general) from that given, for instance, by Theorem 2.7 from \cite{AgKu}, 

\begin{Pro}\label{uniqueness} Let $T\in\mathcal{B}(\H)$ be a real normal operator, and let $G_\R: Bor_{E_\C}(\mathfrak{S})\break\mapsto\mathcal{B}(\H)$ be a real spectral 
measure with the following properties.

(a) $T G_\R(A)\H\subseteq G_\R(A)\H$ for all $A\in Bor_{E_\C}(\mathfrak{S})$.

(b) We have the inclusion $\sigma_\C(T\vert G_\R(A)\H)\subseteq A$ for all $A=[A]\subseteq\mathfrak{S}$.

Then the real spectral measure $G_\R$ coincides with the real spectral measure $E_\R$, which is  the restriction to the $\sigma$-algebra 
$Bor_{E_\C}(\mathfrak{S})$  of  the spectral measure $E_\C$ of the complex extension $T_\C$ of $T$. 
\end{Pro}

{\it Proof}\, For every $A\in Bor_{E_\C}(\mathfrak{S})$ we set $G_\C(A)=G_\R(A)_\C$.  Denoting by $\mathcal{F}^c(\mathfrak{S})$ the family of all 
closed subset in $\rm{Bor}_{E_\C}(\mathfrak{S})$, we get  what is called a {\it pseudoring of closed subsets} of   $\mathfrak{S}$ (see \cite{Vas0}, Definition IV.1.1). Moreover, the restriction 
$$
\mathcal{F}^c(\mathfrak{S})\ni A\mapsto G_\C(A)\H_\C\subset\H_\C,
$$
is a {\it spectral capacities} attached to the normal operator $T_\C$ (see  \cite{Vas0}, Definitions IV.1.5 and IV.1.6), because a normal operator is {\it decomposable} (see
\cite{CoFo} for the original definition, which is a special case of Definition IV.1.6 
from   \cite{Vas0}). Then a particular case of Theorem  IV.1.9 from \cite{Vas0}
implies the equality  $ G_\C(A)=E_\C(A)$ for all 
$A\in\mathcal{F}^c(\mathfrak{S})$. The regularity of the measures leads to the equality $ G_\C(A)=E_\C(A)$ for all $A\in Bor_{E_\C}(\mathfrak{S})$, whence $ G_\R(A)=E_\R(A)$ for all $A\in Bor_{E_\C}(\mathfrak{S})$. 

\noindent{\bf QED}


\begin{Rem}\rm
 A real normal operator may be said to be $\mathcal{F}^c(\mathfrak{S})$-{\it decomposable}, in the spirit of Definition IV.1.6 from \cite{Vas0}, while the original concept of {\it decomposable operator} was introduced in \cite{Foi} (see also \cite{CoFo}). We also recall that the concept of {\it spectral capacity} was introduced in \cite{Apo}.
\end{Rem}

 \subsection{Unbounded Normal Operators in Real Hilbert Spaces}

As in the bounded case, some important  properties of a real normal unbounded operator may and will be obtained from the well-known propertis of its complex normal extension.
For some properties of unbounded normal operators used in what follows, we refer to \cite{Rud} (especially the subsections 13.22-13.25).
For some notation and comments see also the Subsection 2.1.

If the operator $T$, possibly unbounded, is a normal operator in the real Hilbert space  $\H$, then its complex extension $T_\C$ is a normal operator in $\H_\C$ and has a unique spectral measure 
$E_\C: \rm{Bor}(\mathfrak{S})\mapsto\mathcal{B}(\H_\C)$, whose values are orthogonal projections commuting with $T_\C$, that is $E_\C(A) T_\C \xi= T_\C E_\C(A)\xi$ for all $\xi\in D(T_\C )$ and $A\in\rm{Bor}(\mathfrak{S})$, where 
$\mathfrak{S}=\sigma(T_\C)$. We also denote by $\{\mu_{\xi,\eta};\xi,\eta\in\H_\C\}$ the family of complex  measures, associated  with the spectral measure $E_\C$.

With this notation, given an arbitrary function $f\in\mathfrak{B(S)}$, we put
\begin{equation}\label{funcdom}
\mathcal{D}_f=\{\xi\in \H_\C; \int_\mathfrak{S}\vert f\vert^2 
d\mu_{\xi,\xi}<\infty\}.
\end{equation}

\begin{Rem}\label{genunbop}\rm

Citing  Theorems 13.24 and  13.33 from \cite{Rud}, we recall some important properties, to be used in the following.   

 Let $T: D(T)\subset\H\mapsto\H$ be a possibly unbounded normal operator, let $T_\C:D(T_\C)\subseteq\H_\C\mapsto\H_\C$ be its complex extension, and let $E_\C:{\rm Bor}(\mathfrak{S})\mapsto \mathcal{B}(\H_\C)$ be the  spectral measure of $T_\C$, where  $\mathfrak{S}=\sigma(T_\C)$. Then we have
$$
\langle T_\C\xi,\eta\rangle_{\H_\C}=\int_\mathfrak{S} z d\mu_{\xi,\eta}(z),\,\xi\in D(T_\C),\eta\in\H_\C,
$$
where $\mu_{\xi,\eta}$ is the measure $\langle E_{\C}(.)\xi,\eta\rangle_{\H_\C}.$

Moreover, for every $S\in\B(\H_\C)$ such that $ST_\C \xi=T_\C S\xi$ for all $\xi\in D(T_\C)$ we also have $E_\C(A)S=SE_\C(A)$ for all Borel sets $A\subset\mathfrak{S}$.

In addition, we have a functional calculus with unbounded Borel functions for the unbounded normal operator $T_\C$,  given by
\begin{equation}\label{int_ubdfc}
\langle\Psi(f)\xi,\eta\rangle_{\H_\C}=\int_\mathfrak{S}f d\mu_{\xi,\eta}, \, f\in\mathfrak{B(S)},\,\,\xi\in \mathcal{D}_f,\,\,\eta\in\H_\C
\end{equation}
\end{Rem}

The next result is a version of Lemma \ref{inv_cond}, valid for unbounded operators.


\begin{Lem}\label{inv_cond_ubd} With $T,\,T_\C$ and $E_\C$ as above,  we have the equality $E_\C(A)^\flat=E_\C(A)$ for some $A\in{\rm Bor}(\mathfrak{S})$ if and only if $\chi_A=\chi_{A^c}$ 
in $L^\infty(\mathfrak{S},E_\C)$.
\end{Lem}

{\it Proof}\, We follow the lines of the proof of Lemma \ref{inv_cond}, so the map $E_\C^\flat:{\rm Bor}(\mathfrak{S})\mapsto\mathcal{B}(\H_\C)$, with $E_\C^\flat(A)=E_\C(A)^\flat$, is a 
spectral measure. A similar argument used in the quoted proof leads to the equality
$$
\int_\mathfrak{S} z d\mu^\flat_{\xi,\eta}=\langle  T_\C^*\xi,\eta\rangle_{\H_\C}
\,\,\,\forall\,\, \xi,\eta\in D(T_\C)=D(T_\C^*),
$$ 
where, as before,  $\mu^\flat_{\xi,\eta}$
is the scalar measure  $\langle E_\C^\flat(*)\xi,\eta\rangle_{\H_\C}$, 
showing that  $E_\C^\flat$ is precisely the spectral measure of the normal operator
$T^*_\C$. 

As in the proof of Lemma  \ref{inv_cond}, the spectral measure of the adjoint  $T^*_\C$ of the normal operator $T_\C$
can also be obtained by regarding   $T^*_\C$ as a function of $T_\C$, via the map
$\C\ni z\mapsto\bar{z}\in\C$, and leading to the equality  $E_\C^\flat(A)= E_\C(A^c)$ for all $A\in{\rm Bor}(\mathfrak{S})$. Consequently, 
$E_\C(A^c)=CE_\C(A)C$.  

The rest of the proof regards only the spectral measures $E_\C$ and $E_\C^\flat$,
and can be performed as in the proof of Lemma \ref{inv_cond}, leading to the desired conclusion.

\noindent{\bf QED}

As in the bounded case, we set 
$$\rm{Bor}_{E_\C}(\mathfrak{S})=\{A\in{Bor}(\mathfrak{S});E_\C(A^c)=E_\C(A)\}=\{A\in{Bor}(\mathfrak{S}); A=A^c\},
$$
which is a $\sigma$-algebra.

 A functional calculus  associated to a real unbounded normal operator with functions from the space
$\mathfrak{B}_{s,E_\C}(\mathfrak{S})$ (see Remark \ref{varprop}(3)) is given by the following result. See also \cite{AgKu1} for a different approach.

\begin{Thm}\label{ubdfc} Let $T: D(T)\subseteq\H\mapsto\H$ be a real normal operator, let $T_\C:D(T_\C)\subseteq\H_\C\mapsto\H_\C$ be its complex extension, and let $E_\C:{\rm Bor}(\mathfrak{S})\mapsto \mathcal{B}(\H_\C)$ be the spectral measure of $T_\C$, where  $\mathfrak{S}=\sigma(T_\C)$. Let also $\{\mu_{\xi,\eta},\xi,\eta\in\H_C\}$ be the family of complex measures associated to the spectral measure $E_\C$.

Then, setting $\mathcal{D}_{f,\H}=\mathcal{D}_f\cap\H$, 
for every  function $f\in\mathfrak{B}_{s,E_\C}(\mathfrak{S})$ we define
the map $\Phi(f)$ via the equality
\begin{equation}\label{phimap}
\langle\Phi(f)x,y\rangle_\H=\int_\mathfrak{S}fd\mu_{x,y},\,\,x\in\mathcal{D}_{f,\H},\, 
y\in\H,
\end{equation}
which is a normal operator with domain $D(\Phi(f))= \mathcal{D}_{f,\H}$,
satisfying 
$$
\Vert \Phi(f)x\Vert^2_\H=\int_\mathfrak{S}\vert f\vert^2d\mu_{x,x},\,\, x\in \mathcal{D}_{f,\H}.
$$ 
Moreover,
$ \Phi(f)\Phi(g)\subseteq \Phi(fg)$, with $D( \Phi(f)\Phi(g))=
\mathcal{D}_{g,\H}\cap\mathcal{D}_{fg,\H}$. When  $D( \Phi(f)\Phi(g))= \mathcal{D}_{g,\H}$,
we actually have $ \Phi(f)\Phi(g)=\Phi(fg)$.
\end{Thm}

{\it Proof}\, Let $f\in\mathfrak{B}_{s,E_\C}(\mathfrak{S})$ be bounded. 
Then the operator $\Psi(f)$, given by the equation (\ref{int_ubdfc}) is a bounded normal
operator, by Theorem 13.24 from \cite{Rud}. Because $f$ is bounded, it is possible to approximate it with elementary functions, as in Remark \ref{varprop}(4), leading to
the conclusion that the space $\H$ is invariant under the operator $\Psi(f)$, via the 
corresponding argument from the proof of  Theorem \ref{bdfc0}. Moreover, the formula
$$
\langle\Psi(f)x,y\rangle_\H=\int_\mathfrak{S}fd\mu_{x,y},\,\, x,y\in\H
$$ 
clearly holds.

If $f\in\mathfrak{B}_{s,E_\C}(\mathfrak{S})$ is not bounded, we define the 
sets $\mathfrak{S}_n=\{z\in\mathfrak{S};\vert f(z)\vert\le n\}$, for every integer
$n\ge1$, which are conjugate symmetric. Setting $f_n=\chi_{\mathfrak{S}_n}f$,
we have
$$
\lim_{n\to\infty}\Vert\Psi(f_n)x-\Psi(f)x\Vert^2_\H\le\lim_{n\to\infty}\int_\mathfrak{S}\vert f_n-f\vert^2 d\mu_{x,x}=0,\,\, x\in \mathcal{D}_{f,\H},
$$
by Legesque's theorem of dominated convergence (as in formula (7) from the proof
of Theorem 13.24 in \cite{Rud}). Consequently, $\Psi(f)x\in\H$.

We therefore set
$\Phi(f)=\Psi(f)\vert \mathcal{D}_{f,\H}$ for all  $f\in\mathfrak{B}_{s,E_\C}(\mathfrak{S})$, which is precisely the map from (\ref{phimap}). Note that we have
$$
\Vert \Phi(f)x\Vert^2_\H=\Vert \Psi(f)x\Vert^2_\H=
\int_\mathfrak{S}\vert f\vert^2d\mu_{x,x},\,\, x\in \mathcal{D}_{f,\H},
$$ 
via the corresponding property of the spectral measure $E_\C$, proved in Theorem 13.24 from \cite{Rud}.

Next, if $x\in\mathcal{D}_{g,\H}\cap\mathcal{D}_{fg,\H}$, then $ \Phi(f)\Phi(g)x= 
\Psi(f)\Psi(g)x=\Psi(fg)x=\Phi(fg)x$, because $\mathcal{D}_{g,\H}\cap\mathcal{D}_{fg,\H}\subset\mathcal{D}_{g}\cap\mathcal{D}_{fg}$. Then, if $\mathcal{D}_{g,\H}\cap\mathcal{D}_{fg,\H}=\mathcal{D}_{g,\H}$, we must have  
$ \Phi(f)\Phi(g)=\Phi(fg)$, again as in in Theorem 13.24 from \cite{Rud}.

\noindent{\bf QED}

  As in he case of Theorem \ref{bdfc0}, Theorem \ref{ubdfc} associates  every real unbounded normal operator 
$T\in\mathcal{C}(\H)$ with a spectral measure $E_\R$, which, according to Lemma \ref{inv_cond_ubd}, should be defined on the $\sigma$-algabra 
$Bor_{E_\C}(\mathfrak{S})$ (see Remark 6). In fact, we have the following.


\begin{Cor} Let $T: D(T)\subseteq\H\mapsto\H$ be a real normal operator, let $T_\C:D(T_\C)\subseteq\H_\C\mapsto\H_\C$ be its complex extension, and let $E_\C:{\rm Bor}(\mathfrak{S})\mapsto \mathcal{B}(\H_\C)$ be the spectral measure of $T_\C$, where  $\mathfrak{S}=\sigma(T_\C)$. Let also $\{\mu_{\xi,\eta},\xi,\eta\in\H_C\}$ be the family of scalar measures associated to the spectral measure $E_\C$.

Then the restriction $E_\R$ of  the complex spectral measure $E_\C$ to $\rm{Bor}_{E_\C}(\mathfrak{S})$ is a real spectral measure 
with values in $\mathcal{B}(\H)$,  such that 
$$
\langle Tx,y\rangle_\H=\int_\mathfrak{S} zd\mu_{x,y},\,\, x\in D(T), y\in\H,
$$
and $\mu_{x,y}(A)=\langle E_\R(A)x,y\rangle_\H$ for all $A\in\rm{Bor}_{E_\C}(\mathfrak{S})$.  

In addition, if $S\in\mathcal{B}(\H)$ and $STx=TSx$ for all $x\in D(T)$, then $SE_\R(A)=E_\R(A)S$ for all  $A\in\rm{Bor}_{E_\C}(\mathfrak{S})$.
\end{Cor}

{\it Proof}\, Most of the assertions follow from Theorem \ref{ubdfc}, and we leave the details to the reader. Let us only prove the last one. 
If  $S\in\mathcal{B}(\H)$ and $STx=TSx$ for all $x\in D(T)$, we must have $S_\C T_\C \xi=
T_\C S_\C \xi$ for all $\xi\in D(T_\C)$. Therefore, $S_\C E_\C(A)=
E_\C(A) S_\C $ for all $A\in{\rm Bor}(\mathfrak{S})$, via Theorem 13.33 from \cite{Rud}.  Passing to restrictions to $\H$,
we obtain $SE_\R(A)=E_\R(A)S$ for all  $A\in\rm{Bor}_{E_\C}(\mathfrak{S})$

\noindent{\bf QED}

A uniqueness result of the real spectral measure attached to a real normal operator
is also valid for unbounded operators.


\begin{Pro}\label{unb_uniqueness} Let $T\in\mathcal{C}(\H)$ be a real normal operator, and let $G_\R:\rm{Bor}_{E_\C}(\mathfrak{S})\break\mapsto\mathcal{B}(\H)$ be a real spectral 
measure with the following properties.

(a)  The inclusion $T(G_\R(A)\H\cap D(T))\subseteq G_\R(A)\H$ holds for all $A\in\rm{Bor}_{E_\C}(\mathfrak{S})$.

(b) The inclusion $\sigma_\C(T\vert G_\R(A)\H\cap D(T))\subseteq A$ holds  for every $A\subseteq\mathfrak{S}$ compact.

Then the real spectral measure $G_\R$ coincides with the real spectral measure obtained 
by the restriction to the $\sigma$-algebra $\rm{Bor}_{E_\C}(\mathfrak{S})$  of  the spectral measure $E_\C$  of the complex extension $T_\C$ of $T$. 
\end{Pro}

{\it Proof}. Fixing $A\in \rm{Bor}_{E_\C}(\mathfrak{S})$  a compact set, as we have the inclusion
$\sigma_\C(T\vert G_\R(A)\H\cap D(T))\subseteq A$, we must have  that   $G_\R(A)\H\subseteq D(T)$ and that $T_A=T\vert  G_\R(A)\H$ is bounded. Applying Proposition \ref{uniqueness}
to the bounded normal operator $T_A$, we deduce that its spectral measure 
$G_{\R,A}(B)=G_\R(A\cap B), B\in \rm{Bor}_{E_\C}(\mathfrak{S})$ is uniquely 
determined. Then the regularity of the spectral measure $G_\R$ leads to its uniqueness.

\noindent{\bf QED}

\section{Normal Operators in Quaternionic\\ Hilbert Spaces} 

\subsection{Quaternionic Algebra as a Hilbert Space}

First of all let us recall some  known definitions and elementary
facts (see for instance \cite{CoSaSt}, Section 4.6, and/or
\cite{Vas1}). 

Let $\Hh$ be the abstract algebra of quaternions, which is the four-dimensional $\R$-algebra with 
unit $1$, generated by the ''imaginary units`` $\{\bf{j,k,l}\}$,  
which satisfy
$$
{\bf jk=-kj=l,\,kl=-lk=j,\,lj=-jl=k,\,jj=kk=ll}=-1.
$$

We may assume that $\Hh\supset\R$ identifying every  number
$x\in\R$ with the element $x1\in\Hh$.

The algebra $\Hh$ has a natural involution
$$
\Hh\ni{\bf x} = x_0+x_1{\bf j}+x_2{\bf k}+x_3{\bf l}\mapsto
{\bf x}^*= x_0-x_1{\bf j}-x_2{\bf k}-x_3{\bf l}\in\H.
$$

For an arbitrary quaternion ${\x}= x_0+x_1{\bf j}+x_2{\bf k}+x_3{\bf l},\,\,x_0,x_1,x_2,x_3\in\R$, we set $\Re{\bf x}=x_0=
({\bf x}+{\bf x}^*)/2$, and $\Im{\bf x}=x_1{\bf j}+x_2{\bf k}+x_3{\bf l}=({\bf x}-{\bf x}^*)/2$, that is, the {\it real} and 
{\it imaginary part} of ${\bf x}$, respectively.

The algebra $\Hh$ may be regarded as a real Hilbert space, with the inner product
$$
\langle \x, \y \rangle_\Hh=\sum_{m=0}^3 x_m y_m,\,\, \x= x_0+x_1{\bf j}+x_2{\bf k}+x_3{\bf l},\,
 \y= y_0+y_1{\bf j}+y_2{\bf k}+y_3{\bf l}\in\Hh,
$$
denoted simply by $\langle \x, \y \rangle$ when no confusion is possible.

A direct computation shows that
$$
\langle \q \x, \y \rangle=\langle \x, \q^*\y \rangle,
\,\,\langle \x \q, \y \rangle=\langle\x,\y\q^*\rangle\,\, \x,\y, \q\in\Hh.
$$

The natural norm associated to the inner product $\langle \x, \y \rangle_\Hh$ is given by
$$
\Vert {\bf x}\Vert_\Hh=\sqrt{x_0^2+x_1^2+x_2^2+x_0^2},\,\,{\bf x}= x_0+x_1{\bf j}+x_2{\bf k}+x_3{\bf l},\,\,x_0,x_1,x_2,x_3\in\R,
$$
which is a multiplicative norm, simply denoted by $\Vert {\bf x}\Vert$.

Note that ${\bf x}{\bf x}^*={\bf x}^*{\bf x}=\Vert{\bf x}\Vert^2$, implying, in particular, that every element ${\bf x}\in\Hh\setminus\{0\}$ is invertible, and ${\bf x}^{-1}=
\Vert {\bf x}\Vert^{-2}{\bf x}^*$.
\medskip

As an element of a real algebra, each ${\bf x}\in\Hh$ has a complex spectrum given by $\sigma_\C({\bf x})=
\{\Re{\bf x}\pm i\Vert\Im{\bf x}\Vert\}$ (see \cite{Vas1} for details). In particular, if $\lambda\in\sigma_\C({\bf x})$, then $\vert\lambda\vert=\Vert {\bf x}\Vert$.
\smallskip

We now consider the complexification $\C\otimes_\R\Hh$
of the $\R$-algebra $\Hh$ (see also \cite{GhMoPe}), which will be
identified with the  direct sum $\M=\Hh+i\Hh$. 
Of course, the algebra $\M$ contains the complex field $\C$.  Moreover, in the algebra $\M$, the elements of $\Hh$ commute with all complex numbers.  In particular, the ''imaginary units`` 
$\bf j,k,l$ of the algebra $\Hh$  are  independent of and commute with the  imaginary unit $i$ of the complex plane $\C$.

The algebra $\M$ is an involutive one, whose involution is given by 
$(\x+i\y)^*=\x^*-i\y^*$ for all $\x, \y \in\Hh$.

In the algebra $\M$, there also exists a natural conjugation defined
by $\bar{\bf a}={\bf b}-i{\bf c}$, where ${\bf a}={\bf b}+i{\bf c}$ is arbitrary in $\M$, with ${\bf b},{\bf c}\in\Hh$ (see also
\cite{GhMoPe}). Note that $\overline{\bf a+b}=\bar{\bf a}+\bar{\bf b}$, and  $\overline{\bf ab}=\bar{\bf a}\bar{\bf b}$, in particular  $\overline{r\bf a}=r\bar{\bf a}$ for all ${\bf a},{\bf b}\in\M$, and $r\in\R$.  Moreover, $\bar{{\bf a}}={\bf a}$ if and only if ${\bf a}\in\Hh$, which is a useful characterization of the elements of $\Hh$ among those of $\M$. 

The algebra $\M$ may also be regarded as a (complex) Hilbert space, with the inner product given by

$$
\langle \x+i\y, {\bf u}+i{\bf v}\rangle_\M=\langle \x, {\bf u}\rangle+\langle \y, {\bf v}\rangle +i\langle \y, 
{\bf u}\rangle-i\langle \x, {\bf v}\rangle
$$
for all $\x,\y,{\bf u},{\bf v}\in\Hh$, which is the natural extension of the inner oroduct of $\Hh$. It will be 
simply denoted by $\langle \x+i\y, {\bf u}+i{\bf v}\rangle$, when no confusion is possiblde. 

It is also clear that
$$
\langle \q(\x+i\y), {\bf u}+i{\bf v}\rangle_\M=\langle \x+i\y, \q^*({\bf u}+i{\bf v})\rangle_\M,
$$
$$
\langle (\x+i\y)\q, {\bf u}+i{\bf v}\rangle_\M=\langle \x+i\y, ({\bf u}+i{\bf v})\q^*\rangle_\M,
$$
for all $\x,\y,\q,{\bf u},{\bf v}\in\Hh$

\subsection{Quaternionic Hilbert spaces}

\begin{Rem}\label{Hspace}\rm
Following \cite{CoSaSt}, a {\it right $\Hh$-vector space} 
$\mathcal{V}$ is a real vector space 
having a right multiplication with the elements of $\Hh$, that is, 
$x1=x,\,(x+y){\bf q}=x{\bf q}+y{\bf q},\,x({\bf q}+{\bf s})=
x{\bf q}+x{\bf s},\, x({\bf q}{\bf s})=(x{\bf q}){\bf s}$
for all $x,y\in\mathcal{V}$ and ${\bf q},{\bf s}\in\Hh$.

In a similar way, one defines the concept of a {\it left $\Hh$-vector space}. A real vector space $\mathcal{V}$ will be said to be an {\it $\Hh$-vector space} if it is simultaneously a right $\Hh$- and a left $\Hh$-vector space.

If $\H$ is a real Hilbert space, which is a right $\Hh$-vector space, the operator 
$T\in\mathcal{B}(\H)$ is {\it right $\Hh$-linear} if 
$T(x{\bf q})=T(x){\bf q}$ for all $x\in\H$ and
$\q\in\Hh$. The set of all right $\Hh$-linear operators will be 
denoted by $\mathcal{B}^{\rm r}(\mathcal{\H})$, which is, in particular, a unital real algebra. Nevertheless,as noticed in \cite{CoSaSt}, it is the framework of  $\Hh$-vector spaces an appropriate one for the study of right $\Hh$-linear operators.

If $\H$ is $\Hh$-vector space which is also a Hilbert space with the inner product
$\langle *,*\rangle$ and norm $\Vert*\Vert$, assuming that
\begin{equation}\label{HH-space}
\langle \q  x,y\rangle=\langle x,\q^*y\rangle,\,\langle x\q,y\rangle=\langle x,y\q^*\rangle,\,\, x,y\in\H,\,\q\in\Hh,
\end{equation}
then $\H$ is said to be a {\it Hilbert $\Hh$-space}. In this case, we have that $ R_\q\in \mathcal{B}(\H)$, and 
the map $\Hh\ni\q\mapsto R_\q\in 
\mathcal{B}(\H)$ is norm continuous, where $R_{\bf q}$ is the right multiplication of the elements of $\H$ by a given quaternion ${\q}\in\Hh$. In fact, this map is actually an isometry because we have 
$\Vert  R_\q x\Vert=\Vert\q\Vert\Vert x\Vert$ for all $\q\in\Hh$ and $x\in\H$.

Similarly, if  $L_{\q}$ is the left multiplication of the elements of 
$\H$ by the quaternion ${\q}\in\Hh$, we have
 $ L_\q\in \mathcal{B}(\H)$ for all ${\bf q}\in\Hh$, because $\Vert  L_\q x\Vert=\Vert\q\Vert\Vert x\Vert$ for all $\q\in\Hh$ and $x\in\H$.
Note also that 
$$
\mathcal{B}^{\rm r}(\H)=\{T\in\mathcal{B}(\H);TR_\q=R_\q T,\,\q\in\Hh\}.
$$

In particular, $L_\q\in\mathcal{B}^{\rm r}(\H)$ for all $\q\in\Hh$.   
\end{Rem}

\begin{Rem}\label{rems}\rm (1) Note that the quaternionic algebra $\Hh$ is a simple
example of a Hilbert $\Hh$-space because it is a Hilbert space and the
equation (\ref{HH-space}) is satisfied in $\Hh$, as seen in the previous subsection.

(2) Let $\H$ be a Hilbert $\Hh$-space and let us
 consider the complexification $\H_\C$ of $\H$. Because $\H$ is an $\Hh$-bimodule, the space $\H_\C$ is actually an $\M$-bimodule, via the multiplications
$$
(\q+i\s)(x+iy)=\q x-\s y+i(\q y+\s x),
(x+iy)(\q+i\s)=x\q-y\s+i(y\q+x\s),
$$ 
for all $\q+i\s\in\M,\,\q,\s\in\Hh,\,x+iy\in\H_\C,\,x,y\in\H$. Moreover, if the operator $T\in\mathcal{B}(\H)$ is  right $\Hh$-linear, the operator 
$T_\C$ is right $\M$-linear, that is $T_\C((x+iy)(\q+i\s))=
T_\C(x+iy)(\q+i\s)$ for all $\q+i\s\in\M,\,x+iy\in\H_\C$, via a direct computation.

(3) Let $C$ be the natural  conjugation of $\H_\C$. As in the real case, for every
$S\in \mathcal{B}(\H_\C)$, we put $S^\flat=CSC$. The left and  right
multiplication with the quaternion $\q$ on $\H_\C$ will be also denoted by $L_\q,R_\q$, respectively, as elements of $\mathcal{B}(\H_\C)$. We set
$$
\mathcal{B}^{\rm r}(\H_\C)=\{S\in \mathcal{B}(\H_\C); SR_\q=R_\q S,\,\q\in\Hh\},
$$
which is a unital complex algebra containing all operators $L_\q,\q\in\Hh$.
Note that if $S\in\mathcal{B}^{\rm r}(\H_\C)$, then $S^\flat\in\mathcal{B}^{\rm r}(\H_\C)$. Indeed, because $CR_\q=R_\q C$, we also have  $S^\flat R_\q=R_\q S^\flat$. In fact, as we have
$(S+S^\flat)(\H)\subset\H$ and $i(S-S^\flat)(\H)\subset\H$, it folows that  
the algebras $\mathcal{B}^{\rm r}(\H_\C),\,\mathcal{B}^{\rm r}(\H)_\C$
are isomorphic, and they will be often identified, where 
$\mathcal{B^{\rm r}(\H)}_\C=\mathcal{B^{\rm r}(\H)}+i\mathcal{B^{\rm r}(\H)}$ is the complexification of $\mathcal{B}^{\rm r}(\H)$, which is also a unital complex Banach algebra.
\end{Rem}

\subsection{Bounded Normal Operators in Quaternionic Hilbert Spaces}

Let $\H$ be a Hilbert $\Hh$-space, and let $\H_\C$ be its complexification. Being, in particular, a real Hilbert space, the adjoint $T^*\in\mathcal{B}(\H)$ of an operator 
$T\in\mathcal{B}(\H)$ and  its normality are defined as in the framework of 
real operators. We are mainly interested by normal operators 
$T\in\mathcal{B}^{\rm r}(\H)$. In this case we also have $T^*\in\mathcal{B}^{\rm r}(\H)$ and, of course, $T$ is normal if $TT^*=T^*T$. 

If $T\in\mathcal{B}^{\rm r}(\H)$ is normal, then $T_\C\in\mathcal{B}^{\rm r}(\H_\C)$
is normal, and it has a spectral measure $E_\C$ defined on ${\rm Bor}(\mathfrak{S})$,
having a  priori values in $\mathcal{B}(\H_\C)$, where
$ \mathfrak{S}=\sigma(T_\C)$. Because the operators $R_\q (\q\in\Hh)$
commute with $T_\C$, they should also commute with its spectral measure,
implying that the spectral measure $E_\C$ takes values actually in $\mathcal{B}^{\rm r}(\H_\C)$.

Then we have the following representation theorem for bounded quaternionic normal operators.


\begin{Thm}\label{Qbdfc} Let $T\in\mathcal{B}^{\rm r}(\H)$ be a normal operator, let $T_\C\in \mathcal{B}^{\rm r}(\H_\C)$ be its complex extension, and let $E_\C:{\rm Bor}(\mathfrak{S})
\mapsto \mathcal{B}^{\rm r}(\H_\C)$ be the spectral measure of $T_\C$, where $ \mathfrak{S}=\sigma(T_\C)$.

 We consider the map $\Phi:L^\infty_s (\mathfrak{S},E_\C)
\mapsto\mathcal{B}(\H)$,
given by the formula
$$
\langle\Phi(f)x, y \rangle=\int_\mathfrak{S} fd\mu_{x,y},\,\, f\in L^\infty_s (\mathfrak{S},E_\C),\, x,y\in\H,
$$
where $\{\mu_{x,y};x,y\in\H\}$ is the corresponding subset of the family of   measures 
$\{\mu_{\xi,\eta}; \xi,\eta\in\H_\C\}$ associated with the spectral measure $E_\C$. The map $\Phi$ is a unital real  algebra isomorphic  morphism from $ L^\infty_s (\mathfrak{S},E_\C)$ into 
$\mathcal{B}^{\rm r}(\H)$.  Moreover,

(i) For every polynomial with real coeficients $p(z,\bar{z}),\,z\in\C,$ one has $\Phi(p)=
p(T_\C,T_\C^*)$.

(ii) $\Phi(f)$ is a real normal operator and $\Phi(\bar{f})=\Phi(f)^*$ for all 
$f\in L^\infty_s (\mathfrak{S},E_\C)$.

(iii) $\Vert \Phi(f)x\Vert^2=\int_\mathfrak{S}\vert f\vert^2d\mu_{x,x}$ for all 
$f\in L^\infty_s (\mathfrak{S},E_\C)$ and $x\in\H$. 

(iv) $\sigma(\Phi(f))\subseteq [f(\mathfrak{S})]$ for all 
$f\in L^\infty_s (\mathfrak{S},E_\C)$.

(v) If  $S\in\mathcal{B}^{\rm r}(\H)$ and $ST=TS$, then $S\Phi(f)=\Phi(f)S$ 
 for all $f\in L^\infty_s (\mathfrak{S},E_\C)$.
\end{Thm}

{\it Proof}\, This a version of Theorem \ref{bdfc0}, whose proof is similar to the quoted result, via some minor modifications. For instance, 
when $f=\Sigma_{j\in J}(r_j\chi_{A_j}+is_j\theta_{A_j})$ is an elementary function,  we must have
$$
f(T_\C)R_\q=  \Sigma_{j\in J}(r_jE_\C(A_j)+is_j\theta_{A_j}(T_\C))R_\q=
R_\q f(T_\C),\,\, \q\in\Hh,
$$
because the spectral measure $E_\C$ is  $ \mathcal{B}^{\rm r}(\H_\C)$-valued. 
This commutation property also holds by passing to limits, showing that 
$f(T_\C)\in  \mathcal{B}^{\rm r}(\H_\C)$ for all $f\in L^\infty_s (\mathfrak{S},E_\C)$.
Consequently, $\Phi(f)=f(T_\C)\vert\H\in\mathcal{B}^{\rm r}(\H)$ for all $f\in L^\infty_s (\mathfrak{S},E_\C)$. 

Other details are left to the reader.

\noindent{\bf QED}


\begin{Rem}\rm Let $T\in\mathcal{B}^{\rm r}(\H)$ be a normal operator, let $T_\C\in \mathcal{B}^{\rm r}(\H_\C)$ be its complex extension, and let $E_\C:{\rm Bor}(\mathfrak{S})
\mapsto \mathcal{B}^{\rm r}(\H_\C)$ be the spectral measure of $T_\C$, where $ \mathfrak{S}=\sigma(T_\C)$. As in Corolary \ref{bdfc1}, we have the following. 

The restriction $E_\R$ of the complex spectral measure $E_\C$ to the $\sigma$-algebra $\rm{Bor}_{E_\C}(\mathfrak{S})$ is a
real spectral measure with values in $\mathcal{B}^{\rm r}(\H)$, 
with
 $$\langle E_\R(A)x,y\rangle=\int_\mathfrak{S}\chi_A d\mu_{x,y}, \,x,y\in\H,$$ and 
 
$$
 \langle Tx,y\rangle=\int_\mathfrak{S} z d\mu_{x,y},\,\,x,y\in\H.
 $$ 
where the family of measures 
$\{\mu_{x,y}:x,y\in\H\}$ on $Bor(\mathfrak{S})$, may be regarded
as linear functionals on the real Banach space  $C_s(\mathfrak{S})$. 
\end{Rem}

We also have a uniqueness result corresponding to Proposition \ref{uniqueness}, having a similar proof, which will be omitted.  

\begin{Pro}\label{Quniqueness} Let $T\in\mathcal{B}^{\rm r}(\H)$ be a real normal operator, and let $G_\R:\rm{Bor}_{E_\C}(\mathfrak{S})\break\mapsto\mathcal{B}^{\rm r}(\H)$ be a real spectral 
measure with the following properties.

(a) $T G_\R(A)\H\subseteq G_\R(A)\H$ for all $A\in\rm{Bor}_{E_\C}(\mathfrak{S})$.

(b) We have the inclusion $\sigma_\C(T\vert G_\R(A)\H)\subseteq A$ for all $A=[A]\subseteq\mathfrak{S}$.

Then the real spectral measure $G_\R$ coincides with the real spectral measure $E_\R$  obtained 
by the restriction to the $\sigma$-algebra $\rm{Bor}_{E_\C}(\mathfrak{S})$  of  the spectral measure $E_\C$  of the complex extension $T_\C$ of $T$. 
\end{Pro}

\subsection{Unbounded Normal Operators in Quaternionic\\ Hilbert Spaces}

Let $\H$ be a  Hilbert $\Hh$-space, which is, in particular, a real Hilbert space. 
Let also $T$ be an operator from the family 
$\mathcal{C}(\H)$. Therefore,  the adjoint $T^*\in\mathcal{C}(\H)$ of $T$ and its normality  are defined as in the framework of  real operators. 

The operator  $T\in\mathcal{C}(\H)$ is said to be {\it right $\Hh$-linear} if  $R_\q D(T)\subset
D(T)$, and $T(x{\bf q})=T(x){\bf q}$ for all $x\in D(T)$ and
$\q\in\Hh$. The family  of right $\Hh$-linear operators from  $\mathcal{C}(\H)$  will be 
denoted by $\mathcal{C}^{\rm r}(\mathcal{\H})$. A direct calculation shows that, if $T\in \mathcal{C}^{\rm r}(\mathcal{\H})$, we have $\langle(T^*R_\q)x,y\rangle=\langle(R_\q T^*)x,y\rangle$ for all $x\in D(T^*),y\in D(T)$, showing that $T^*\in \mathcal{C}^{\rm r}(\mathcal{\H})$.

As before, we are particularly interested by normal operators 
$T\in\mathcal{C}^{\rm r}(\H)$. Then $T^*\in\mathcal{C}^{\rm r}(\H)$, and $T$
is normal if $D(T^*)=D(T)$, $D(TT^*)=D(T^*T)$, and $TT^*=T^*T$. 

We now consider  the complexification $\H_\C$ of $\H$. If $T\in\mathcal{C}(\H)$, then
$T_\C\in\mathcal{C}(\H_\C)$, and $D(T_\C)=D(T)_\C$. Moreover, if $T\in \mathcal{C}^{\rm r}(\mathcal{\H})$, then $T_\C\in\mathcal{C}^{\rm r}(\mathcal{\H_\C})$, where the latter space consists of those operators from 
$\mathcal{C}(\H_\C)$ which are right $\Hh$-linear.

When  the operator $T\in\mathcal{C}(\H)$ is  normal, and so is its complex extension $T_\C\in\mathcal{C}(\H_\C)$ is also normal, having a unique spectral measure 
$E_\C: \rm{Bor}(\mathfrak{S})\mapsto\mathcal{B}(\H_\C)$, whose values are orthogonal projections commuting with $T_\C$, as in Subsection 2.4. Let also
$\{\mu_{\xi,\eta};\xi,\eta\in\H_\C\}$ be the family of complex valued associated measures with the spectral
measure $E_\C$. If  $T\in \mathcal{C}^{\rm r}(\mathcal{\H})$, and so $T_\C\in\mathcal{C}^{\rm r}(\mathcal{\H_\C})$   because $R_\q T_\C \xi=T_\C R_\q \xi$ for all $\xi\in D(T_\C)$ and $\q\in\Hh$,  we must actually  have 
$E_\C: \rm{Bor}(\mathfrak{S})\mapsto\mathcal{B}^{\rm r}(\H_\C)$, as follows from
Theorem 13.33 from \cite{Rud} or Proposition 5.26 from\cite{Sch}.
\medskip

The next result is a functional calculus with a large class of Borel functions, valid for unbounded
normal operators in Hilbert $\Hh$-spaces. The notation $\mathcal{D}_f$ is given by the equality (\ref{funcdom}).


\begin{Thm}\label{Qubdfc} Let $\H$ be an  $\Hh$-space, and let 
$T\in\mathcal{C}^{\rm r}(\H)$ be a normal operator. Let also $T_\C\in\mathcal{C}^{\rm r}(\H_\C)$ be its complex extension, and let $E_\C:{\rm Bor}(\mathfrak{S})\mapsto \mathcal{B}^{\rm r}(\H_\C)$ be the spectral measure of $T_\C$, where  $\mathfrak{S}=\sigma(T_\C)$. 
 
 Setting $\mathcal{D}_{f,\H}=\mathcal{D}_f\cap\H$, 
for every function $f\in\mathfrak{B}_{s,E_\C}(\mathfrak{S})$ we put
$$
\langle\Phi(f)x,y\rangle=\int_\mathfrak{S}fd\mu_{x,y},\,\, x\in \mathcal{D}_{f,\H}, y\in\H,
$$ 
which is a normal operator in $\mathcal{C}^{\rm r}(\H)$,  with domain $D(\Phi(f))= \mathcal{D}_{f,\H}$,
satisfying 
$$
\Vert \Phi(f)x\Vert^2=\int_\mathfrak{S}\vert f\vert^2d\mu_{x,x},\,\, x\in \mathcal{D}_{f,\H},
$$
where $\{\mu_{x,y};x,y\in\H\}$ is the corresponding subset of the family of   measures 
$\{\mu_{\xi,\eta}; \xi,\eta\in\H_\C\}$ associated with the spectral measure $E_\C$. 

Moreover,
$ \Phi(f)\Phi(g)\subseteq \Phi(fg)$, with $D( \Phi(f)\Phi(g))=
\mathcal{D}_{g,\H}\cap\mathcal{D}_{fg,\H}$. When  $D( \Phi(f)\Phi(g))= \mathcal{D}_{g,\H}$,
we actually have $ \Phi(f)\Phi(g)=\Phi(fg)$.

In addition, if $S\in\mathcal{B}^{\rm r}(\H)$ and $STx=TSx$ for all $x\in D(T)$, then $SE_\R(A)=E_\R(A)S$ for all  $A\in\rm{Bor}_{E_\C}(\mathfrak{S})$. 
\end{Thm}

{\it Proof}\, The proof of this assertion is similar to that of Theorem \ref{ubdfc}, except for 
some details. For instance, it is clear that for every $E_\C$-stem function $f\in\mathfrak{B}_{s,E_\C}(\mathfrak{S})$ the operator $\Phi(f)$ is normal as an
element of  $\mathcal{C}(\H)$. In fact, it actually belongs to $\mathcal{C}^{\rm r}(\H)$.
Of course, this is true when $f$ is bounded, as in the proof of Theorem \ref{Qbdfc}.
For an arbitrary unbounded function $f\in\mathfrak{B}_{s,E_\C}(\mathfrak{S})$, we refine an argument from the proof of Theorem \ref{ubdfc}.

We set  $\mathfrak{S}_n=\{z\in\mathfrak{S};\vert f(z)\vert\le n\}$, for every integer
$n\ge1$, which are conjugate symmetric. Setting also $f_n=\chi_{\mathfrak{S}_n}f$,
as we have $\lim_{n\to\infty}\Psi(f_n) x=\Psi(f) x$, and 
 $R_\q\Psi(f_n) x=\Psi(f_n)R_\q x$, it follows
$$
\lim_{n\to\infty}R_\q\Psi(f_n) x=\lim_{n\to\infty}\Psi(f_n)R_\q x=R_q \Psi(f) x,
$$
showing that $R_\q x\in D(\Psi(f))$, and $R_\q\Psi(f)x=\Psi(f) R_\q x$ for all $x\in \mathcal{D}_f\cap\in\H$. Therefore, $\Psi(f)\in\mathcal{C}^{\rm r}(\H) $. 

Other details are left to the reader.

\noindent{\bf QED}

\begin{Rem}\rm Let $\H$ be an  $\Hh$-space, and let 
$T\in\mathcal{C}^{\rm r}(\H)$ be a normal operator. Let also $T_\C\in\mathcal{C}^{\rm r}(\H_\C)$ be its complex extension, and let $E_\C:{\rm Bor}(\mathfrak{S})\mapsto \mathcal{B}^{\rm r}(\H_\C)$ be the spectral measure of $T_\C$, where  $\mathfrak{S}=\sigma(T_\C)$.

Then the restriction $E_\R$ of  the complex spectral measure $E_\C$ to $\rm{Bor}_{E_\C}(\mathfrak{S})$ is a real spectral measure $E_\R$ 
with values in $\mathcal{B}^{\rm r}(\H)$,  such that 
$$
\langle Tx,y\rangle=\int_\mathfrak{S} zd\mu_{x,y},\,\, x\in D(T), y\in\H,
$$
where $\mu_{x,y}(A)=\langle E_\R(A)x,y\rangle$ for all $A\in\rm{Bor}_{E_\C}(\mathfrak{S})$.

\end{Rem} 

As in some previous cases, we still have a uniqueness result. 

\begin{Pro}\label{Qubd_uniqueness}  Let $\H$ be  $\Hh$-space, let $T\in\mathcal{C}^{\rm r}(\H)$ be a real normal operator, and let $G_\R:\rm{Bor}_{E_\C}(\mathfrak{S})\mapsto\mathcal{B}^{\rm r}(\H)$ be a real spectral 
measure with the following properties.

(a)  The inclusion $T(G_\R(A)\H\cap D(T))\subset G_\R(A)\H$ holds for all $A\in\rm{Bor}_{E_\C}(\mathfrak{S})$.

(b) The inclusion $\sigma_\C(T\vert G_\R(A)\H\cap D(T))\subset A$ holds  for every $A\subseteq\mathfrak{S}$ compact.

Then the real spectral measure $G_\R$ coincides with the real spectral measure $E_R$  obtained 
by the restriction to the $\sigma$-algebra $\rm{Bor}_{E_\C}(\mathfrak{S})$  of  the spectral measure $E_\C$  of the complex extension $T_\C$ of $T$. 
\end{Pro}

 \section{Other Examples} 


\begin{Exa}\rm One of a simplest examples of a normal operator of a
Hilbert $\Hh$-space is the left multiplication operator $L_\q$ acting on the Hilbert $\Hh$-space $\Hh$, given by $L_\q\s=\q\s$ for
all $\s\in\Hh$, where $\q\in\Hh$ is a fixed quaternion.  The complex extension of this operator to the Hilbert space $\Hh_\C=\M$ is given
by a similar formula and we have $\sigma_\C(L_\q)= 
\sigma_\C({\q})=
\{\Re{\q}\pm i\Vert\Im{\q}\Vert\}$, as recalled in Subsection 3.1. The adjoint $L_\q^*$ of $L_\q$ is clearly the 
multiplication operator $L_{\q^*}$, and so the operator $L_\q$ on 
the Hilbert space $\Hh$ is normal. Let us exhibit the spectral representation of the normal operator $L_\q$. 
We give some arguments
in this sense but most of the details can be found in \cite{Vas1}, Subsections 3 and 4.  

If ${\bf q}\in\H$ and  ${\Im \bf q}\neq 0$, setting 
$\mathfrak{s}_{\tilde{\bf q}}=\tilde{\bf q}\Vert\tilde{\bf q}\Vert
^{-1 }$, where $\tilde{\bf q}=\Im\bf q$, we define 
$\iota_\pm(\mathfrak{s}_{\tilde{\bf q}})=
2^{-1}(1\mp i\mathfrak{\mathfrak{s}_{\tilde{\bf q}}})$ in $\M$. 

We put $s_\pm(\q)=\Re{\q}\pm i\Vert\Im{\q}\Vert$, which are the eigenvalues of $L_\q$, and denote
by $P_\pm(\q)$ the corresponding spectral projections, which are given by
$$
P_\pm({\bf q}){\bf a}=\iota_\pm(\mathfrak{s}_{\tilde{\bf q}}){\bf a},\,\,{\bf a}\in\M.
$$
In fact, we have  $P_+({\bf q})P_-({\bf q})=P_-({\bf q})P_+({\bf q})=0$, and 
$P_+({\bf q})+P_-({\bf q})$ is the identity on $\M$. Moreover,
$P_\pm(\q)^*=P_\pm(\q)$, as one can easily see. This shows that 
the family
$\{0,P_+(\q),P_-(\q), I_\M\}$ defines a spectral measure on the finite set $\sigma_\C(L_\q)$. The associated functional calculus for $L_\q$ is defined by the formula
$$
f(L_\q)=f(s_+(\q))P_+(\q)+f(s_-(\q))P\_(\q)\in\M,
$$
where $f:\sigma_\C(L_\q)\mapsto\C$. This map is linear and 
multiplicative, and it is $\Hh$-valued if and 
only if $f(s_+(\q))=\overline{f(s_-(\q))}$.

Identifying $L_\q$ with $\q$, for a given  ${\bf q}\in\H$ with 
${\Im \bf q}\neq 0$ we have, in the algebra $\M$, the spectral decomposition
$$
\q=s_+(\q)\iota_+(\mathfrak{s}_{\tilde{\bf q}})+
s_-(\q)\iota_-(\mathfrak{s}_{\tilde{\bf q}}),
$$
where $\iota_\pm(\mathfrak{s}_{\tilde{\bf q}})$ are
self-adjoint commuting idempotents, whose product is null, and their sum is the  identity.

When ${\bf q}\in\R$, the corresponding spectral 
decomposition is trivial, equal to the identity of 
$\M$. 

Because $s_\pm(\q^*)=s_\pm(\q)$ and 
$\iota_\pm(\mathfrak{s}_{\tilde{\bf q}^*})=
\iota_\mp(\mathfrak{s}_{\tilde{\bf q}})$, we have

$$
\q^*=s_-(\q)\iota_+(\mathfrak{s}_{\tilde{\bf q}})+
s_+(\q)\iota_-(\mathfrak{s}_{\tilde{\bf q}}). 
$$

Therefore, for a polynomial $p(z,\bar{z})$, we infer that
$$
p(\q,\q^*)=p(s_+(\q),s_-(\q))\iota_+(\mathfrak{s}_{\tilde{\bf q}})+p(s_-(\q),s_+(\q))\iota_-(\mathfrak{s}_{\tilde{\bf q}}),
$$
provided $\Im(\q)\neq 0$. 

\end{Exa}


\begin{Exa}\rm
We shall give a basic example of a quaternionic normal operator, defined by the multiplication with the 
independent quaternionic variable in a Hilbert space of 
quaternionic-valued square integrable functions. 

First of all, we recall some definitions and properties from  \cite{Vas1} (see especially Remark 4), also adding a few simple observations. A subset $S$ in $\Hh$ is said to be {\it spectrally saturated} if whenever for an element 
$\q\in\Hh$ we have $\sigma_\C(\q)=\sigma_\C(\s)$ for some $\s\in S$, then $\q\in S$.
Setting $\mathfrak{S}(\Omega)=\cup_{\q\in\Omega}\sigma_\C(\q)$, for an arbitrary subset $\Omega\subseteq\Hh$, we have 
that $\Omega$ is spectrally saturared if an only if  $\mathfrak{S}(\Omega)$ is conjugate symmetric.

Let $\K\subset\Hh$ be a spectrally saturated compact subset. Then the conjugate symmetric set  $\mathfrak{S}(\K)$ is also compact because it is bounded and closed.

Let $\lambda$ be a positive measure on $\K$. Let also $L^2(\K,\Hh)$ be  the space (of equivalence classes) of square integrable functions $f:\K\mapsto\Hh$, regarded as a real vector space. Using the inner product defined on $\Hh$ (see Subsection 3.1), we consider  the scalar product
$$
\langle f,g\rangle=\int_\K\langle f(\q),g(\q)\rangle_\Hh d\lambda(\q)\,\,\,\, f,g\in L^2(\K,\Hh),
$$
on the vector space $L^2(\K,\Hh)$, which becomes a Hilbert $\Hh$-space, simply denoted  by $\H$. Note that $\H_\C$ is isomorphic to the complex Hilbert space $L^2(\K,\M)$, whose inner product is induced by that of $\M$. In other words,

$$
\langle \phi,\psi\rangle_{\H_\C}=\int_\K\langle \phi(\q),\psi(\q)\rangle_\M d\lambda(\q)\,\,\,\, \phi,\psi\in L^2(\K,\M).
$$

Let $T$ be the map given by $Tf(\q)=\q f(\q)$ for all $f\in\H$ and $\q\in\K$, that is, the left multiplication with the quaternionic independent variable $\q$ in the space $\H$. Then 
$T\in\mathcal{B}^{\rm r}(\H)$
and it is a normal operator, because $T^*f(\q)=\q^*f(\q)$ for all $f\in\H$ and $\q\in\K$, and obviously 
$TT^*=T^*T$. 

Let us compute the complex spectrum $\sigma_\C(T)=\sigma(T_\C)$ of $T$, where the complex extension $T_\C$ of $T$
is clearly the left multiplication with the independent quaternionic variable $\q$ in the space $\H_\C$.

We shall prove that the spectrum  $\sigma(T_\C)$ is given by the conjugate symmetric set 
$\mathfrak{S}(\K)$.  We must show that a complex number $\zeta\notin\sigma(T_\C)$ if and only if 
$\zeta\notin\mathfrak{S}(\K)$ for all $\q\in\K$.  

If $\zeta\notin\sigma(T_\C)$, then 
$(\zeta-\q)[(\zeta-T_\C)^{-1}\iota](\q)=1$ for all 
$\q\in\K$, where $\iota(\q)=1$.  In other 
words, $\zeta -\q$ has a right inverse in $\M$. 
Because we also have $\bar{\zeta}\notin\sigma(T_\C)$, applying 
this equality to $\bar{\zeta}$ and $\q^*$, we obtain 
$(\bar{\zeta}-\q^*)[(\bar{\zeta}-T_\C)^{-1}\iota](\q^*)=1$, implying 
that $\{[(\bar{\zeta}-T_\C)^{-1}1](\q^*)\}^*$ is a left inverse in $\M$ of $\zeta-\q$. Therefore, $\zeta\notin\mathfrak{S}(\K)$  for all 
$\q\in\K$.

If $\zeta\notin\mathfrak{S}(\K)$, the function 
$\q\mapsto(\zeta-\q)^{-1}$ is well defined on $\K$ and the left multiplication with this function on $\H_\C$ is equal to the operator 
$(\zeta-T_\C)^{-1}$, via the fact that the operator $\zeta-T_\C$ is 
bijective.  Hence $\zeta\notin\sigma_\C(T)$. Writing explicitly the formula of the spectrum of an arbitrary quaternion
we obtain $\mathfrak{S}(\K)=\{\Re(\q)\pm i\Vert\Im(\q)\Vert;\q\in\K\}=\sigma(T_\C)$. 

Let us start a discussion about the spectral structure of the operator $T_\C$.

If $p(z,\bar{z})$ is a polynomial with complex coefficients restricted to $\mathfrak{S}(\K)$, the action of the operator $p(T_\C,T_\C^*)$ in
$\H_\C$ is the multiplication with the function $p(\q,\q^*),\,\q\in\K$.
Because the polynomials $p(z,\bar{z})$ with complex coefficients are uniformly dense in $C(\mathfrak{S}(\K))$, for a function $f\in C(\mathfrak{S}(\K))$, the action of the operator $f(T_\C)$ should be given
by the multiplication with the function $f(\q)$. Using the notation from the previous example, we have
$$
f(\q)= \lim_{n\to\infty}[p_n(s_+(\q),s_-(\q))\iota_+(\mathfrak{s}_{\tilde{\bf q}})+p_n(s_-(\q),s_+(\q))\iota_-(\mathfrak{s}_{\tilde{\bf q}})],
$$
for each $\q\in\K, \Im(\q)\neq 0$, where $(p_n)_n$ is a sequence 
of polynomials uniformly convergent to $f$. Of course, if  
$\Im(\q)=0$, $f(\q)$ is a constant.

As a matter of fact, we may define 
an operator $f(T_\C)$ for a bounded Borel function as a multiplication operator with the function $f(\q),\q\in\K$. 
To define this operator we may use the poinwise definition 
from the previous example, globally extended to $\K$.
Specifically, we ''extend`` the function $f$ to $\K$ by putting
\begin{equation}\label{ext_f}
f(\q)=f(s_+(\q))\iota_+(\mathfrak{s}_{\tilde{\bf q}})+
f(s_-(\q))\iota_-(\mathfrak{s}_{\tilde{\bf q}}), \q\in\K,
\end{equation}
where $f\in L^\infty(\K,\C)$. This $\M$-valued function is 
measurable because the functions $\K\ni\q\mapsto s_\pm(\q)\in\C$
are continuous, the functions $\K\setminus\R\ni\q\mapsto  
\iota_\pm(\mathfrak{s}_{\tilde{\bf q}})\in\M$ are also continuous and the $\C$-valued function $f$ is measurable.
Moreover, $\Vert f(\q)\vert\le 2\Vert f\Vert_\infty$ for all
$\q\in\K$. Therefore, we have a map from $L^\infty(\K,\C)$ into $B^r(\H_\C)$ associating
each  complex-valued function $f$ with an operator denoted by $f(T_\C)$, which is the multiplication operator with the 
$\M$-valued function 
$f$ given by (\ref{ext_f}), that is, $(f(T_\C)\phi)(\q)=f(\q)\phi(\q)$ for all points $\q\in\K$, and for 
 an arbitrary function $\phi\in\H_\C$. Using the properties obtained in the previous example, we deduce that this map is unital, linear and multiplicative. In particular, if 
 $A\in Bor( \mathfrak{S}(\K))$, we may consider the operator 
 $\chi_A(T_\C)$ as a multiplication with the function 
 $\chi_A(\q),\q\in\K$. Note that 
 $$
 \chi_A(\q)=\iota_+(\mathfrak{s}_{\tilde{\bf q}})\,\,{\rm if}\,\,s_+(\q))\in A,\,s_-(\q))\notin A;
 $$
 $$
 \chi_A(\q)=\iota_-(\mathfrak{s}_{\tilde{\bf q}})\,\,{\rm if}\,\,s_-(\q))\in A,\,s_+(\q))\notin A;
 $$
 $$
 \chi_A(\q)=1\,\,{\rm if}\,\, \sigma_\C(\q)\subseteq A; 
  \chi_A(\q)=0\,\,{\rm if}\,\, \sigma_\C(\q)\cap A=\emptyset.
 $$
 This shows that the map $\q\mapsto\chi_A(\q)$ is a self-adjoint
 idempotent in $B^{\rm r}(\H_\C)$.   Setting
$E(A)=\chi_A(T_\C)$, we obtain an operator-valued measure $E:Bor(\mathfrak{S}(\K))\mapsto B^{\rm r}(\H_\C)$, which will be shown to the spectral measure of $T_\C$. 

We put $E(\emptyset)=0$ and we clearly
have $E(\K)$ to be the identity on $\H_\C$. If $A,B\in Bor(\mathfrak{S}(\K))$ are arbitrary, the multiplicativity of the map $\q\mapsto f(\q)$ shows that $E(A\cap B)=E(A)E(B)$. In addition, if 
$(A_k)_{k\ge1}$ is a countable family of mutually disjoint subsets from $Bor(\mathfrak{S}(\K))$, then 

$$\langle\sum_{k\ge1}E(A_k)\phi,\psi\rangle_{\H_\C}=
\langle\bigcup_{k\ge1}E(A_k)\phi,\psi\rangle_{\H_\C},\,\,
\phi,\psi\in\H_\C,
$$ 
via the pointwise convergence of the sequence 
$(\sum_{k\ge1}^m E(A_k))_{m\ge1}$ to $(\cup_k E(A_k))$ and the Lebesgue theorem of dominated convergence. In other words, the 
map $E$ is a spectral measure. To see that it is the spectral measure of $T_\C$, we first note that the space $E(A)\H_\C$ is invariant under this operator
for all $A\in Bor( \mathfrak{S}(\K))$. In addition, the spectrum 
$\sigma(T_\C\vert E(A)\H_\C)$ is a subset $A=[A]$, because the function
$h(z)=(\lambda-z)^{-1}\chi_A(z), z\in\mathfrak{S}(\K),$ is well defined when $\lambda\notin A$, and has the property 
$(\lambda-\q)h(\q)=\chi_A(\q),\,\q\in\K$. Hence the function $h$
provides  an inverse for $\lambda-T_\C\vert E(A)\H_\C$. It follows that $E$ is precisely the spectral measure of $T_\C$, via Proposition \ref{Quniqueness}. In particular, $F(A)=E(A)\vert\H$ is a spectral measure for the operator $T$, assuming $A=A^c\in 
Bor(\mathfrak{S}(\K))$. In addition, the operator $f(T)$ is well defined on $\H$ if $f$ is a bounded Borel function with the property
$f(\bar{z})=\overline{f(z)}$.

An example of a not necessarily bounded quaternionic normal operator
can be obtained supposing $\K$ to be closed, not necessarily compact. We omit the details.
\end{Exa}


\begin{thebibliography}{xx}


\bibitem{AgKu} S. N. Agrawal, S. H. Kulkarni, A spectral theorem for a normal operator on  real Hilbert space, {\it Acta Sci. Math.} 59 (1994),no. 3-4, 441–451.

\bibitem{AgKu1} S. N. Agrawal, S. H. Kulkarni, Spectral theorem for 
unbounded normal operators in a real Hilbert space; {\it J. Analysis}, 12 (2004), 107-114.

\bibitem{AlCoKi} D. Alpay, F. Colombo, and D. P. Kimsey, The spectral theorem for quaternionic unbounded normal operators based on the S-spectrum, {\it J. Math. Phys.} 57 (2016), no. 2, 023503, 27.

\bibitem{Apo} C. Apostol, Spectral decompositions and functional calculus, {\it Rev. Roum. Math. Pures Appl.} 13 (1968), 1481-1528.

\bibitem{BaZa} A. G. Baskakov and A. S. Zagorskii, Spectral Theory of Linear Relations on Real Banach Spaces.{\it Mathematical Notes} (Russian: Matematicheskie Zametki), 2007, Vol. 81, No. 1, pp. 15-27.

\bibitem{CoFo} I. Colojoar\v a and C. Foia\c s,{\it Theory of Generalized Spectral Operators},  Gordon and Breach, New York, 1968.
 
\bibitem{CoGaKi} F. Colombo, J. Gantner, D. P. Kimsey, {\it Spectral
Theory on the S-Spectrum for Quaternionic Operators},
 Birkh\"auser, 2018.
 
 
\bibitem{Con} J.B. Conway, {\it A Course in Operator Theory}, Graduate Studies in Mathematics, Amer. Math. Soc.,
Providence, RI, 2000
 
 \bibitem{CoSaSt} F. Colombo, I. Sabadini and D. C. Struppa,
{\it Noncommutative Functional Calculus: Theory and Applications of
Slice Hyperholomorphic Functions},  Progress in Mathematics, Vol. 28 Birkh\"auser/Springer Basel AG, Basel, 2011. 

\bibitem{DuSc} N. Dunford and J. T. Schwartz, \textit{Linear Operators,
Part II: Spectral Theory.} Interscience Publishers, New York, London, 1963. \textit{Part III: Spectral Operators.} Wiley-Interscience,  New York, London, 1971.

\bibitem{Foi} C. Foia\c s, Spectral maximal spaces and decomposable operators, {\it Arch. Math.} 14 (1963), 341-369. 

\bibitem{Fug} B. Fuglede, A Commutativity Theorem for Normal Operators, {\it PNAS} 36(1), (1950), 35-40


\bibitem{GhMoPe} R. Ghiloni , V. Moretti
and A. Perotti, Continuous slice functional calculus in quaternionic Hilbert spaces. {\it Rev. Math. Phys.} 25 (2013), no.4, 1350006, 83 p.

\bibitem{GhMoPe1}R. Ghiloni, V. Moretti, A. Perotti, Spectral representations of normal operators in quaternionic
Hilbert spaces via intertwining quaternionic PVMs, {\it Rev. Math. Phys.} 29 (2017).

\bibitem{Good} K. R. Goodearl, {\it Notes on Real and Complex $C^*$-algebras}, Shiva Publishing Limited, England, 1982.  

\bibitem{Goo} R. K. Goodrich, The spectral theorem for real Hilbert space, {\it Acta Sci. Math.} (Szeged) 33 (1972), 123–127.

\bibitem{Grz} M. Grzesiak, Real function algebras and their sets of antisymmetry, {\it Glasnik Matematicki}, 24 (1989), 297–304.

\bibitem{Ing} L. Ingelstam,  Real Banach algebras.{\it Ark. Mat.} 5 
(1964), 239–270 (1964).

\bibitem{Kap} I. Kaplansky,  Normed algebras. {\it Duke. Math. J.} 16, 399-418 (1949). 

\bibitem{KulLim} S. H. Kulkarni and B. V. Limaye , {\it Real Function Algebras}, Marcel Dekker, NewYork, 1992.

\bibitem{Ore} M. N. Oreshina, Spectral decomposition of normal operator in real Hilbert space. {\it Ufa Math. J.} (2017) 9:4 87-96.

\bibitem{RaKu} G. Ramesh, P. Santhosh Kumar, Spectral theorem for quaternionic normal operators: Multiplication form {\it Bull. Sci. math.} 159 (2020) 102840 

\bibitem{Rud} W. Rudin, {\it Functional Analysis}, McGraw-Hill Book Company, New York,
 1973.b

\bibitem{Sch} K. Schm\"udgen, {\it Unbounded Self-adjoint Operators on Hilbert Spaces}, Graduate Texts in Mathematics 265, Springer, Dordrecht, 2012.

\bibitem{Tei} O. Teichm\"uller, Operatoren im Wachsschen Raum. {\it J. Reine Angew. Math.} 174 (1936), 73-124.

\bibitem{Vas0}  F. -H. Vasilescu,  {\it Analytic Functional Calculus and Spectral Decompositions},  D. Reidel Publishing Company, Dordrecht, 1982

\bibitem{Vas1} F.-H. Vasilescu,  Quaternionic Regularity via Analytic Functional Calculus. {\it Integr. Equ. Oper. Theory} 92, 18 (2020). https://doi.org/10.1007/s00020-020-2574-7 .

\bibitem{Vas3} F.-H. Vasilescu,  Spectrum and Analytic Functional Calculus for Clifford Operators via Stem Functions. {\it  Concrete Operators}, vol. 8, no. 1, 2021, pp. 90-113. https://doi.org/10.1515/conop-2020-0115

\bibitem{Vas2} F.-H. Vasilescu, Spectrum and Analytic Functional Calculus in Real and Quaternionic Frameworks. {\it  Pure and Applied Functional Analysis} 7:1 (2022), 389-407.
ISSN 2189-3764.

\bibitem{Vis0} K. Viswanath, Operators on real Hilbert spaces, {\it J. Indian Math. Soc.}, 42 (1978),
1–13.

\bibitem{Vis} K. Viswanath, Normal Operators in Quaternionic Hilbert Spaces. {\it Trans. Amer. Math. Soc.} Vol 162 (1971), 337-350.

\bibitem{Won} T. K. Wong, On Real Normal Operators, {\it ProQuest LLC}, Ann Arbor, MI, 1969, Thesis (Ph.D.)-Indiana University.


\end{thebibliography}
\end{document}